\documentclass[preprint,1p,times]{elsarticle}

\usepackage[utf8]{inputenc}
\usepackage{amsmath}
\usepackage{dutchcal}
\usepackage{amsfonts}
\usepackage{amssymb}
\usepackage{amsthm}
\usepackage[dvipsnames,table]{xcolor}
\usepackage{graphicx}
\usepackage{tikz}
\usetikzlibrary{spy}

\usepackage{booktabs}

\usepackage{pgfplots}
\usetikzlibrary{pgfplots.groupplots}
\usepackage{subcaption}

\usepackage{bm}
\usepackage{mathtools}
\usepackage{enumerate}
\usepackage{hyperref}
\hypersetup{
    colorlinks = true,
}
\usepackage[capitalize]{cleveref}

\usepackage{algorithm}
\usepackage{algpseudocode}

\usepackage{physics}

\newcommand{\pp}{{\textbf{r}}} 
\newcommand{\surf}{{\widehat{\textbf{r}}}} 
\newcommand{\breg}{{\textbf{s}}} 
\newcommand{\reg}{{s}} 
\newcommand{\bfun}{{r}} 

\definecolor{col1}{rgb}{0.83529411764705885,0.24313725490196081,0.30980392156862752}
\definecolor{col2}{rgb}{0.9882352941176471,0.55294117647058827,0.34901960784313718}
\definecolor{col3}{rgb}{0.99607843137254903,0.8784313725490196,0.54509803921568623}
\definecolor{col4}{rgb}{0.90196078431372551,0.96078431372549022,0.59607843137254901}
\definecolor{col5}{rgb}{0.59999999999999998,0.83529411764705885,0.58039215686274515}
\definecolor{col6}{rgb}{0.19607843137254899,0.53333333333333333,0.74117647058823533}

\usepackage{lineno}

\biboptions{sort&compress}

\pgfplotscreateplotcyclelist{mark list mod}{
every mark/.append style={solid,fill=\pgfplotsmarklistfill,mark size=1pt,fill opacity=0.5},mark=*\\
every mark/.append style={solid,fill=\pgfplotsmarklistfill,mark size=1pt,fill opacity=0.5},mark=square*\\
every mark/.append style={solid,fill=\pgfplotsmarklistfill,mark size=1pt,fill opacity=0.5},mark=triangle*\\
every mark/.append style={solid,fill=\pgfplotsmarklistfill,mark size=1pt,fill opacity=0.5},mark=diamond*\\
every mark/.append style={solid,fill=\pgfplotsmarklistfill,mark size=1pt,fill opacity=0.5},mark=pentagon*\\
}

\pgfplotsset{
    legend pos = south east,
    colormap/jet,
    cycle multiindex* list={[samples of colormap={5}]\nextlist mark list mod},
    width = \linewidth,
    height=0.25\textheight,
    grid=major,
    legend cell align={left},
    ylabel near ticks,
    xlabel near ticks,
    legend style={fill=white, fill opacity=0.6, draw opacity=1,text opacity=1, font=\small},
    enlargelimits
}

\theoremstyle{definition}

\crefname{rem}{Remark}{Remarks}

\crefname{ex}{Example}{Examples}

\begin{document}
\begin{frontmatter}

\title{Isogeometric analysis for multi-patch structured Kirchhoff-Love shells}

\author[label1]{Andrea Farahat\corref{cor2}}
\ead{andrea.farahat@ricam.oeaw.ac.at}
\author[label2,label3]{Hugo M. Verhelst\corref{cor2}}
\ead{h.m.verhelst@tudelft.nl}
\author[label4]{Josef Kiendl}
\ead{josef.kiendl@unibw.de}
\author[label5]{Mario Kapl\corref{cor1}}
\ead{m.kapl@fh-kaernten.at}
\cortext[cor1]{Corresponding Author. }
\cortext[cor2]{These authors have made an equal contribution to this paper. The breakdown of the contributions is given in \cref{sec:credit}.}

\address[label1]{Johann Radon Institute for Computational and Applied Mathematics, Austrian Academy of Sciences, Linz, Austria}
\address[label2]{Department of Applied Mathematics, Delft University of Technology, The Netherlands}
\address[label3]{Department of Maritime and Transport Technology, Delft University of Technology, The Netherlands}
\address[label4]{Institute of Engineering Mechanics Structural Analysis, Universit\"at der Bundeswehr M\"unchen, Munich, Germany}
\address[label5]{ADMiRE Research Center - Additive Manufacturing, Intelligent Robotics, Sensors and Engineering, School of Engineering and IT, Carinthia University of Applied Sciences, Villach, Austria}

\begin{abstract}
We present an isogeometric method for Kirchhoff-Love shell analysis of shell structures with geometries composed of multiple patches and which possibly possess extraordinary vertices, i.e. vertices with a valency different to four. The proposed isogeometric shell discretisation is based on the one hand on the approximation of the mid-surface by a particular class of multi-patch surfaces, called analysis-suitable~$G^1$~\cite{CoSaTa16}, and on the other hand on the use of the globally $C^1$-smooth isogeometric multi-patch spline space~\cite{FaJuKaTa22}. We use our developed technique within an isogeometric Kirchhoff-Love shell formulation~\cite{kiendl-bletzinger-linhard-09} to study linear and non-linear shell problems on multi-patch structures. Thereby, the numerical results show the great potential of our method for efficient shell analysis of geometrically complex multi-patch structures which cannot be modeled without the use of extraordinary vertices.
\end{abstract}


\end{frontmatter}



\section{Introduction}

Isogeometric analysis~\cite{Hughes2005} (IGA) is a powerful numerical framework for solving partial differential equations by employing the same (rational) spline function space for representing the geometry and the solution space of the considered partial differential equation. In the last years, it has been widely used as an effective tool for the analysis of shells in various engineering disciplines including aerospace, automotive, maritime, civil and biomedical engineering. Thereby, the high continuity of the employed splines allows to achieve a high accuracy and efficiency in the numerical simulation. 
Shell analysis is a field where IGA had an especially high impact, and many innovative shell formulations have been developed, including Kirchhoff-Love shells \cite{kiendl-bletzinger-linhard-09,Alaydin2021}, Reissner-Mindlin shells \cite{benson2011large,Zou2020,Oesterle2016,Dornisch2013,Dornisch2014,Dornisch2016}, hierarchic shells \cite{Echter2013,Oesterle2017}, and solid shells \cite{Bouclier2013,Hosseini2013,caseiro_on_2014,caseiro_assumed_2015,Leonetti2018}. 

In this work, we will focus on the use of the Kirchhoff-Love shell element~\cite{kiendl-bletzinger-linhard-09}. It is a popular element due to its rotation-free nature, which is enabled by the use of second-order derivatives of the underlying basis functions, and which implies a low number of element-wise degrees of freedom. The requirement for the existence of second-order derivatives of the underlying basis of the isogeometric Kirchhoff-Love shell elements can be easily satisfied for spline functions within a single patch. But in case of multi-patch geometries with possibly extraordinary vertices, which are typically needed for modeling complex shell structures, the $C^1$ continuity across the patch interfaces has to be additionally enforced. There exist two main strategies for imposing $C^1$-smoothness across the interfaces of multi-patch structured Kirchhoff-Love shells.

The first approach couples the neighboring patches in a weak sense, which means that the used isogeometric discretisation space is just approximately $C^1$-smooth instead of exactly $C^1$-smooth. The so-called bending strip method \cite{Kiendl2010} was the first developed technique for weak patch coupling and falls under the category of penalty methods. 
It uses the concept of a stiff strip to connect patches with arbitrary angle and  requires a conforming patch interface. Besides the bending strip method, different penalty methods \cite{Herrema2019,Coradello,Leonetti2020,Zhao2022}, Nitsche coupling techniques \cite{Guo2015,Guo2019,Nguyen-Thanh2017,Nguyen2014,Apostolatos2014,Hu2018}, and other methods which are reviewed in \cite{Schuss2019} have been developed. While penalty methods provides a simple implementation but require the choice of a parameter, the Nitsche coupling techniques are parameter-free but need larger implementation efforts. 
An alternative to penalty and Nitsche methods are mortar methods \cite{Horger2019,Dornisch2015,Bouclier2017}. They are parameter-free and also provide flexibility for the patch connection angle as well as for non-conforming coupling but require the finding of Lagrangian multipliers. Contrary to the penalty, Nitsche and mortar methods, where the variational formulation is altered to establish the weak $C^1$ coupling across patches, the methods~\cite{SaJu21,Takacs2022,Weinmuller2021,WeTa22} directly generate basis functions which are approximately $C^1$. However, the techniques~\cite{SaJu21,Weinmuller2021,WeTa22} are restricted so far to planar domains and to the solving of the biharmonic equation. 

The second strategy for imposing $C^1$-continuity across the patch interfaces couples the neighboring patches in a strong sense, which leads to $C^1$-smooth discretisation spaces that are exactly $C^1$-smooth. Then after constructing a basis for the resulting $C^1$-smooth space, the Kirchhoff-Love shell formulation~\cite{kiendl-bletzinger-linhard-09} can be directly applied. The existing techniques for the Kirchhoff-Love shell analysis can be classified depending on the employed multi-patch parameterisation of the mid-surface of the considered Kirchhoff-Love shell. They are either based on the use of globally $C^1$-smooth multi-patch surfaces with singularities at the extraordinary vertices~\cite{NgPe16,Toshniwal2017}, on the application of globally $C^1$-smooth multi-patch surfaces with specific $G^1$-smooth caps in the vicinity of extraordinary vertices~\cite{NgKaPe15, Pe15-2}, or on the most general case namely on the use multi-patch surfaces which are just $G^1$-smooth\footnote{A surface is $G^1$-smooth if it possesses at each point a well-defined tangent plane, cf. \cite{Pe02}. The $G^1$-continuity of a surface is a weaker condition than the $C^1$-continuity, where at each point the partial derivatives have to be well-defined.}~\cite{ChAnRa18}. For the latter case, the usage of a particular class of $G^1$-smooth multi-patch surfaces, called analysis-suitable-$G^1$ (in short AS-$G^1$)~\cite{CoSaTa16}, allows the construction of $C^1$-smooth multi-patch spline spaces with optimal polynomial reproduction properties. This has been demonstrated so far just for the solving of the biharmonic equation over planar multi-patch parameterisations~\cite{CoSaTa16,KaViJu15,KaBuBeJu16,C1UnstructuredMP,KaSaTa19b} and multi-patch surfaces~\cite{FaJuKaTa22}, and will be extended in this work to the case of multi-patch shell structures. Thereby, the representation of the mid-surface of the shell by an AS-$G^1$ multi-patch surface is not restrictive, since any (approximate) $G^1$-smooth multi-patch surface can be approximated by an AS-$G^1$ multi-patch surface~\cite{FaJuKaTa22,KaSaTa17b}.

In this work, we will develop a novel, simple isogeometric method for the analysis of shell structures composed of several patches with possibly extraordinary vertices. The presented technique will follow the second strategy above, and will rely on three main ingredients, namely first on the representation of the mid-surface of the shell by an AS-$G^1$ multi-patch surface, then on the usage of the globally exactly  $C^1$-smooth multi-patch spline space~\cite{FaJuKaTa22} as discretisation space of the shell and finally on the application of the Kirchhoff-Love shell formulation~\cite{kiendl-bletzinger-linhard-09} for the analysis. An advantage of our proposed approach compared to other strong coupling techniques is that our method can be applied to any $G^1$-smooth multi-patch spline surface, which is then modeled by an AS-$G^1$ multi-patch surface, and is not based on the usage of a particular class of multi-patch surfaces. E.g. the technique~\cite{Toshniwal2017} cannot be applied to smooth multi-patch spline surfaces, which possess boundary vertices, where two patches just meet with $C^0$-smoothness. A benefit of our strong coupling method compared to existing weak coupling techniques is that the employed discretisation space~\cite{FaJuKaTa22} is exactly $C^1$-smooth, which allows to directly apply the Kirchhoff-Love shell formulation~\cite{kiendl-bletzinger-linhard-09}. Hence, it is not necessary e.g. to add penalty terms to the weak form of the problem and to deal with penalty factors as in penalty methods or to use Lagrange multipliers as in mortar methods. Several linear and non-linear benchmark examples will demonstrate the great potential of the presented method for performing analysis of complex Kirchhoff-Love shells.

The outline of the paper is as follows. \cref{sec:KLShell} provides the basics of isogeometric Kirchhoff-Love shell analysis with the focus on the used Kirchhoff-Love shell formulation~\cite{kiendl-bletzinger-linhard-09}. In~\cref{sec:construction}, we introduce the employed $C^1$-smooth multi-patch discretisation, which is based on the one hand on the approximation of the mid-surface of the shell by a particular class of multi-patch surfaces, called AS-$G^1$ multi-patch geometries~\cite{CoSaTa16}, and on the other hand on the application of the $C^1$-smooth spline space~\cite{FaJuKaTa22}. The detailed construction of the used $C^1$-smooth spline space with implementation details is discussed in \ref{app:construction_space_A}.
\cref{sec:numerical_results} presents the novel isogeometric method for the analysis of multi-patch structured Kirchhoff-Love shells
with several linear and non-linear numerical benchmark examples. Finally, we conclude our work in~\cref{sec:conclusion}.

\section{Kirchhoff-Love shell formulation} \label{sec:KLShell}

Based on the works \cite{kiendl-bletzinger-linhard-09,Kiendl2015,Verhelst2021}, we will briefly recall the Kirchhoff-Love shell formulation, which will be used throughout the paper. For the sake of brevity, we will restrict ourselves in this section to a single-patch mid-surface~$\vb{r}:[0,1]^2 \rightarrow \mathbb{R}^3$, but the presented formulation can be simply extended to the employed multi-patch setting introduced in Section~\ref{sec:construction} by just applying it in each case to the single surface patches. Firstly, \cref{subsec:coordinate} will define the shell coordinate system, followed by \cref{subsec:kinematics} where the shell kinematics will be defined accordingly. Finally, in \cref{subsec:variational} the variational form of the shell problem will be given. In the following, we will use Greek indices $\alpha,\beta, \gamma, \delta \in\{1,2\}$ and Latin indices $i,j\in\{1,2,3\}$.

\subsection{Shell coordinate system}\label{subsec:coordinate}

The Kirchhoff-Love shell formulation is defined on the surface $\vb{r}(\xi^1,\xi^2)$ with parametric coordinates $\xi^\alpha$. By the Kirchhoff Hypothesis \cite{Reddy2014}, which implies no shear in the cross section and orthogonality of orthogonal vectors after deformation, any point in the shell $\vb{y}(\xi^1,\xi^2,\xi^3)$ can be represented by a point on the mid-surface and a contribution in normal direction:
\begin{equation}\label{eq:coordinate}
 \vb{y}(\xi^1,\xi^2,\xi^3) = \vb{r}(\xi^1,\xi^2) + \xi^3\vb{a}_3,
\end{equation}
where $\vb{a}_3$ is the unit normal vector to the surface and $\xi^3$ is the through-thickness coordinate. The deformed and undeformed configurations are denoted by $\vb{y},\vb{r}$ and $\mathring{\vb{y}},\mathring{\vb{r}}$, respectively. The covariant basis of the mid-surface and the normal vector are obtained by the partial derivatives of $\vb{r}$ with respect to its parametric coordinates, i.e.
\begin{equation*}\label{eq:conv_basis}
 \vb{a}_\alpha = \pdv{\vb{r}}{\xi^\alpha}, \quad \vb{a}_3 = \frac{\vb{a}_1\times\vb{a}_2}{\abs{\vb{a}_1\times\vb{a}_2}}.
\end{equation*}
In addition, the first and second fundamental forms $a_{\alpha\beta}, b_{\alpha\beta}$  
are defined as follows:
\begin{equation}\label{eq:curvature}
 a_{\alpha\beta} = \vb{a}_\alpha\cdot\vb{a}_\beta,\quad b_{\alpha\beta} = \vb{a}_3\cdot\vb{a}_{\alpha,\beta}=-\vb{a }_{3,\beta}\cdot\vb{a}_\alpha,
\end{equation}
where $\vb{a}_{\alpha,\beta}$ is the Hessian of the surface and $\vb{a}_{3,\alpha}$ is the derivative of the normal vector, which can be obtained by Weingarten's formula $\vb{a}_{3,\alpha}=-b_\alpha^\beta\vb{a}_\beta$ with $b_\alpha^\beta=a^{\alpha\gamma}b_{\gamma\beta}$ as the mixed curvature tensor \cite{Sauer2017}. From \cref{eq:curvature} it can be observed that second-order derivatives are required for evaluation of the second fundamental form. 

\subsection{Shell kinematic and constitutive relations}\label{subsec:kinematics}
The Green-Lagrange strains $E_{\alpha\beta}$ at any point in the shell are defined as: 
\begin{equation}
\begin{aligned}
 E_{\alpha\beta} &= \varepsilon_{\alpha\beta} + \xi_3\kappa_{\alpha\beta},
\end{aligned}
\end{equation}
where $\varepsilon_{\alpha\beta}, \kappa_{\alpha\beta}$ are the membrane strains and curvature change, respectively, which are obtained by the first and second fundamental forms \eqref{eq:curvature} of the undeformed and deformed configurations \cite{kiendl-bletzinger-linhard-09,Kiendl2015}: 
\begin{equation}\label{eq:strains}
\begin{aligned}
 \varepsilon_{\alpha\beta} &= \frac{1}{2}(a_{\alpha\beta} - \mathring{a}_{\alpha\beta}) \\  \kappa_{\alpha\beta} &= b_{\alpha\beta} - \mathring{b}_{\alpha\beta} .
\end{aligned}
\end{equation}
Stresses are represented by the stress resultants $\vb*{n}$ and $\vb*{m}$, corresponding to membrane forces and moments, respectively. Assuming isotropic linear elastic material, they are obtained by
\begin{equation*}
    \begin{aligned}
	    n^{\alpha\beta} &= \mathbb{C}^{\alpha\beta\gamma\delta}:\varepsilon_{\gamma\delta}, \\
	    m^{\alpha\beta} &= \mathbb{C}^{\alpha\beta\gamma\delta}:\kappa_{\gamma\delta}, 
	\end{aligned}
\end{equation*}
where $\mathbb{C}^{\alpha\beta\gamma\delta}$ is the plane stress material tensor. Within this paper, only isotropic linear elastic materials are considered, but the formulation can be easily extended to nonlinear materials as shown in \cite{Kiendl2015}. 

\subsection{Variational formulation}\label{subsec:variational}
The variational formulation for the Kirchhoff-Love shell problem is defined by using the principle of virtual work as in \cite{kiendl-bletzinger-linhard-09,Kiendl2015},
and by following the notations from \cite{Verhelst2021}. Denoting the internal and external energies by $W^{\text{int}}$ and $W^{\text{ext}}$, respectively, the variations of the internal and external work are defined as
\begin{equation*}\label{eq:variational_form}
\begin{aligned}
 \delta W(\vb{u},\delta\vb{u}) &= \delta W^{\text{int}} - \delta W^{\text{ext}}  = \\ & = \int_\Omega \vb*{n}:\delta\bm{\varepsilon} + \vb*{m}:\delta\bm{\kappa} \dd{\Omega} - \int_\Omega \vb{f}\cdot\delta\vb{u} \dd{\Omega}.
\end{aligned}
\end{equation*}
Here, we denote by $\vb{u}$ the displacements, by $\delta\vb{u}$ the virtual displacements, by $\delta\bm{\varepsilon}$ the virtual membrane strain and by $\delta\bm{\kappa}$ the virtual curvature change. \\
Furthermore, $\mathring{\mathbf{r}}$ represents the undeformed mid-surface of the shell, defined on the parametric domain $[0,1]^2$. Note that the metric basis of the deformed configuration is computed on $\mathbf{r}$, which has the same parametric domain as the undeformed configuration. In addition, $\dd{\Omega}=\sqrt{\abs{\mathring{a}_{\alpha\beta}}}\dd\xi^1\dd\xi^2$ is the differential area 
in the undeformed configuration mapped
onto the parametric domain $[0,1]^2$.

Following a Galerkin approach, we represent the shell displacements $\vb{u}$ 
by a finite sum of basis functions 
$\phi_i$ and their coefficients $u_i$,  i.e. $\vb{u} = \sum_{i} u_i \phi_i$. Then, the discrete residual vector $R_i$ is defined by taking the first Gateaux derivative with respect to $u_i$ \cite{kiendl-bletzinger-linhard-09}, i.e.
\begin{equation}\label{eq:residual}
 R_i = F^{\text{int}}_i - F^{\text{ext}}_i = \int_{\Omega} \vb*{n}:\pdv{\bm{\varepsilon}}{u_i}+ \vb*{m}:\pdv{\bm{\kappa}}{u_i} \dd{\Omega} - \int_{\Omega} \vb{f}\cdot\pdv{\vb{u}}{u_i}\dd{\Omega},
\end{equation}
where $\Omega \subset \mathbb{R}^3$ represents the surface domain defined by the mid-surface~$\vb{r}$. Applying a second linearisation of the variational form with respect to $u_j$, the components $K_{ij}$ of the tangential stiffness matrix, can be found as
\begin{equation}\label{eq:jacobian}
\begin{aligned}
 K_{ij} &= K^{\text{int}}_{ij} - K^{\text{ext}}_{ij} = \\
 &= \int_{\Omega} \pdv{\vb*{n}}{u_j}:\pdv{\bm{\varepsilon}}{u_i}+\vb*{n}:\pdv{\bm{\varepsilon}}{u_i}{u_j} + \pdv{\vb*{m}}{u_j}:\pdv{\bm{\kappa}}{u_i} + \vb*{m}:\pdv{\bm{\kappa}}{u_i}{u_j}\dd{\Omega}.
\end{aligned}
\end{equation}

Lastly, for non-linear simulations, Newton iterations are performed for solution $\vb{u}$ and increment $\Delta\vb{u}$ by solving
\begin{equation*}
 K\Delta\vb{u} = -\vb{R}.
\end{equation*}
Throughout these iterations, the deformed geometry $\mathbf{r}$ is updated with the displacement field. 

Due to the appearance of the curvature variations of the curvature tensor $\kappa_{\alpha\beta}$ in \cref{eq:residual,eq:jacobian}, second-order derivatives for the basis functions $\phi_i$ are required for the variational formulation. This means that the basis functions~$\phi_i$ have to be at least globally $C^1$-smooth. Considering a multi-patch structured Kirchhoff-Love shell as in this work, the $C^1$-smoothness can be trivially satisfied within a single surface patch but has to be imposed across the patch interfaces, which gives rise to the globally $C^1$-smooth isogeometric multi-patch spline functions presented in \cref{sec:construction}. 

\section{$C^1$-smooth multi-patch discretisation space} \label{sec:construction}

We will briefly describe the $C^1$-smooth multi-patch space construction~\cite{FaJuKaTa22}, which will be used as discretisation space for the isogeometric Kirchhoff-Love shell analysis of multi-patch geometries with possibly extraordinary vertices in Section~\ref{sec:numerical_results}. Before, we will also present the employed multi-patch structure for the mid-surface~$\pp$ as well as the particular type of multi-patch surface, called AS-$G^1$ multi-patch surface~\cite{CoSaTa16}, which is needed to represent the mid-surface~$\pp$. Details of the construction and implementation of the basis of the $C^1$-smooth spline space will be discussed in \ref{app:construction_space_A}.

\subsection{The multi-patch surface structure}

In the following, let the mid-surface $\pp$ be a $G^1$-smooth conforming multi-patch surface which consists of regular quadrilateral surface patch parameterisations~$\pp^{(i)} \in (\mathcal{S}_h^{\vb{p}, \breg})^3$, $i \in \mathcal{I}_{\Omega}$, given by
 \begin{equation*}
  \pp^{(i)} : [0, 1]^2 \rightarrow \overline{\Omega^{(i)}},
\end{equation*}
where $\mathcal{S}_h^{\vb{p}, \breg}$ is the tensor-product spline space $\mathcal{S}_h^{p, \reg} \otimes \mathcal{S}_h^{p, \reg}$, with $\vb{p}=(p,p)$ and $\breg=(\reg,\reg)$ on the parameter domain $[0,1]^2$, obtained from the univariate spline space $\mathcal{S}_h^{p, \reg}$ of degree $p$, continuity $C^\reg$ and mesh size $h = \frac{1}{k}$ defined on the parameter domain $[0, 1]$. We will require that $p \geq 3, \, 1 \leq \reg \leq p-2$, $k \geq 1$, 
and will further denote by $N_{j}^{p,\reg}$, $j=0,\ldots,n-1$, with $n=p+(k-1)(p-\reg)+1$, the B-splines of the spline space~$\mathcal{S}_h^{p, \reg}$, and by $N^{\vb{p},\breg}_{\vb{j}} = N_{j_1}^{p,\reg}N_{j_2}^{p,\reg}$, $\vb{j}=(j_1,j_2) \in \{0,\ldots,n-1 \}^2$, the associated tensor-product B-splines of $\mathcal{S}_h^{\vb{p}, \breg}$. We will further consider the univariate spline spaces $\mathcal{S}_h^{p, \reg+1}$ and $\mathcal{S}_h^{p-1, \reg}$ with the corresponding B-splines $N_{j}^{p,\reg+1}$, $j=0,\ldots,n_0-1$, $n_0=p+(k-1)(p-\reg-1)+1$ and $N_{j}^{p-1,\reg}$, $j=0,\ldots,n_1-1$, $n_1=p+(k-1)(p-\reg-1)$, respectively.  

\begin{figure}[htbp]
\centering
    \includegraphics[width=0.9\linewidth]{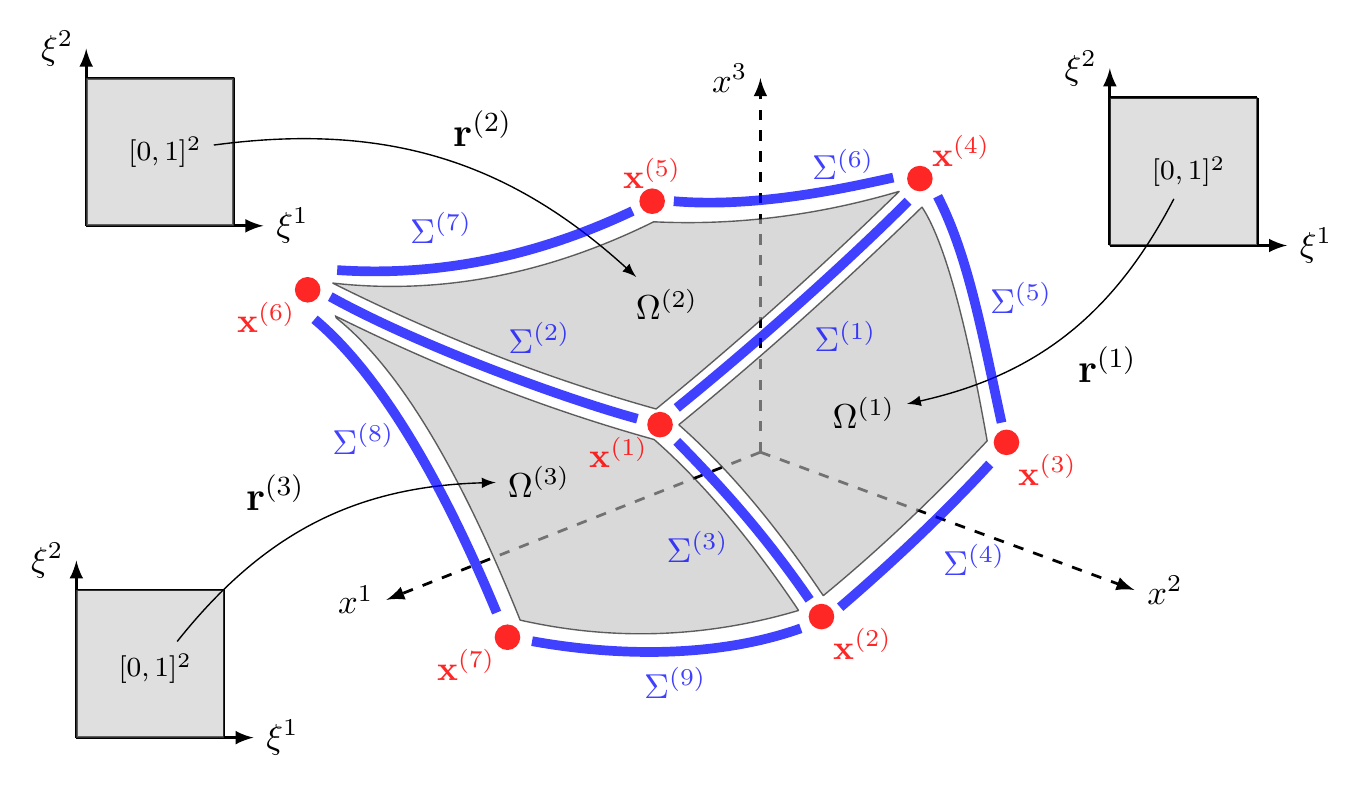}
    \caption{The multi-patch setting demonstrated on the basis of an example of a three-patch surface.}
    \label{fig:geoMap}
\end{figure}

 The multi-patch surface~$\pp$ defines a surface domain~$\Omega \subset \mathbb{R}^3$, which can be represented as the disjoint union of the open quadrilateral surface patches $\Omega^{(i)}$, $i \in \mathcal{I}_{\Omega}$, of open edges (i.e. interface and boundary curves) $\Sigma^{(i)}$, $i \in \mathcal{I}_{\Sigma}$, and of inner and boundary vertices~$\vb{x}^{(i)}, i \in \mathcal{I}_{\chi}$,  
i.e. 
\begin{equation*}
  \Omega = \left( \bigcup_{i \in \mathcal{I}_{\Omega}} \Omega^{(i)}\right) \cup \left( \bigcup_{i \in \mathcal{I}_{\Sigma}} \Sigma^{(i)}\right) \cup \left( \bigcup_{i \in \mathcal{I}_{\chi}} \vb{x}^{(i)}\right),
\end{equation*}
cf. \cref{fig:geoMap}. To distinguish the case of an interface or boundary curve and the case of an inner or boundary vertex, we further divide the index sets~$\mathcal{I}_{\Sigma}$ and $\mathcal{I}_{\chi}$ into $\mathcal{I}_{\Sigma}=\mathcal{I}_{\Sigma}^{\circ} \cup \mathcal{I}_{\Sigma}^{\Gamma}$ and $\mathcal{I}_{\chi}=\mathcal{I}_{\chi}^{\circ} \cup \mathcal{I}_{\chi}^{\Gamma}$, where in both cases the symbol $\Gamma$ denotes the boundary case and the symbol $\circ$ the interface/inner case.

\subsection{AS-$G^1$ multi-patch surfaces} \label{subsec:AS_multi-patch-surface}

In this work, we will model the mid-surface~$\pp$ by a specific $G^1$-smooth multi-patch surface, called analysis-suitable $G^1$ (in short AS-$G^1$) multi-patch surface~\cite{CoSaTa16}. AS-$G^1$ multi-patch surfaces are of great importance, since they are needed to ensure the construction of $C^1$-smooth isogeometric multi-patch spline spaces with optimal polynomial reproduction properties. More precisely, in case that a multi-patch surface is not AS-$G^1$, the approximation power of the resulting $C^1$-smooth spline space over the multi-patch surface can be dramatically reduced, cf.~\cite{KaSaTa17b}. 

We recall first that a $C^0$-smooth multi-patch surface~$\pp$ is $G^1$-smooth if and only if for any two neighboring patches~$\Omega^{(i_1)}$ and $\Omega^{(i_2)}$, $i_1,i_2 \in \mathcal{I}_{\Omega}$, with the common interface curve~$\Sigma^{(i)}\subset \overline{\Omega^{(i_1)}} \cap \overline{\Omega^{(i_2)}}$, there exists functions $\alpha^{(i, i_1)} : [0, 1] \rightarrow \mathbb{R}$, $\alpha^{(i, i_2)} : [0, 1] \rightarrow \mathbb{R}$ and $\beta^{(i)} : [0, 1] \rightarrow \mathbb{R}$ satisfying 
\begin{equation*}
  \alpha^{(i, i_1)}(\xi) \, \alpha^{(i, i_2)}(\xi) > 0, \quad \xi \in [0,1],
\end{equation*} 
and
\begin{equation} \label{ASCOND} 
  \alpha^{(i, i_1)}(\xi) \, \partial_2 \pp^{(i_2)}(\xi, 0) + \alpha^{(i, i_2)}(\xi) \, \partial_1 \pp^{(i_1)}(0, \xi) + \beta^{(i)}(\xi) \, \partial_2 \pp^{(i_1)}(0, \xi) = \vb{0}, \quad \xi \in [0,1],
\end{equation}
cf.~\cite{Pe02}, where $\pp^{(i_1)}$ and $\pp^{(i_2)}$ are the corresponding patch parameterisations of $\Omega^{(i_1)}$ and $\Omega^{(i_2)}$, respectively, both (re)-parameterised (if needed) into the so-called standard form~\cite{FaJuKaTa22}, which means that the common interface curve~$\Sigma^{(i)}$ is given by $\pp^{(i_1)}(0,\xi)=\pp^{(i_2)}(\xi,0)$. (See \ref{appsec:standardform} for details about the standard form.)
Thereby, the functions~$\alpha^{(i,i_1)}$, $\alpha^{(i,i_2)}$ and $\beta^{(i)}$ are uniquely determined up to a common function $\gamma^{(i)}: [0,1] \rightarrow \mathbb{R}$ via the patch parameterisations~$\pp^{(i_1)}$ and $\pp^{(i_2)}$. A $G^1$-smooth multi-patch surface is then called AS-$G^1$ if for each interface curve~$\Sigma^{(i)}$, $i \in \mathcal{I}_{\Sigma}^{(i)}$, the functions~$\alpha^{(i,i_1)}$ and $\alpha^{(i,i_2)}$ can be linear polynomials and if there further exists (non-uniquely determined) linear polynomials $\beta^{(i,i_1)}:[0,1] \rightarrow \mathbb{R}$ and $\beta^{(i,i_2)}:[0,1] \rightarrow \mathbb{R}$ such that 
\begin{equation*}
  \beta^{(i)}(\xi) = \alpha^{(i, i_1)}(\xi) \, \beta^{(i, i_2)}(\xi) + \alpha^{(i, i_2)}(\xi) \, \beta^{(i, i_1)}(\xi) \quad \xi \in [0,1].
\end{equation*}
To uniquely determine the functions $\alpha^{(i, i_1)}$, $\alpha^{(i, i_2)}$, $\beta^{(i, i_1)}$ and $\beta^{(i, i_2)}$, we select those relatively prime functions $\alpha^{(i,i_1)}$ and $\alpha^{(i,i_2)}$ and those functions $\beta^{(i,i_1)}$ and $\beta^{(i,i_2)}$, which minimize
\[
 ||\alpha^{(i,i_1)}-1 ||^{2}_{L_{2}([0,1])} +  ||\alpha^{(i,i_2)} -1 ||^{2}_{L_{2}([0,1])} \quad \mbox{and} \quad  ||\beta^{(i,i_1)} ||^{2}_{L_{2}([0,1])} +  ||\beta^{(i,i_2)} ||^{2}_{L_{2}([0,1])},
\]
respectively, cf.~\cite{C1UnstructuredMP}. In case of a boundary curve~$\Sigma^{(i)} \subset \overline{\Omega^{(i_1)}}$, $i \in \mathcal{I}_{\Sigma}$, we simply set $\alpha^{(i,i_1)}=1$ and $\beta^{(i,i_1)}=0$. 

The design of AS-$G^1$ multi-patch spline surfaces was studied so far for the case of planar domains in~\cite{CoSaTa16,KaBuBeJu16,KaSaTa17b,KaSaTa19b}, and for the case of multi-patch surfaces in~\cite{KaSaTa17b,FaJuKaTa22}. In this work, we will employ the methods~\cite{KaSaTa17b,FaJuKaTa22} to generate AS-$G^1$ multi-patch parameterisations of the mid-surface~$\pp$, e.g. by approximating a non-AS-$G^1$ multi-patch parameterisation of the mid-surface~$\pp$ by an AS-$G^1$ multi-patch geometry. 

\subsection{The specific $C^1$-smooth multi-patch spline space}  \label{subsec:specific_space}

Let the mid-surface~$\pp$ be an AS-$G^1$ multi-patch surface. The space of $C^1$-smooth isogeometric spline functions over $\pp$ is defined as 
\begin{equation*}
\mathcal{V}^1= \{ \phi  \in C^1(\Omega) \; : \; \phi \circ \pp^{(i)} \in \mathcal{S}_{h}^{\vb{p},\breg}, i \in \mathcal{I}_{\Omega} \},
\end{equation*}
or equivalently as
\begin{equation*}
\mathcal{V}^1= \{ \phi  \in L^2(\Omega) \; : \: \phi \circ \pp^{(i)} \in \mathcal{S}_{h}^{\vb{p},\breg}, i \in \mathcal{I}_{\Omega}, \mbox{ and } \bfun_j^{(i,i_1)}(0,\xi) = \bfun_j^{(i,i_2)}(\xi,0), \; \xi \in [0,1] \; j=0,1 , \; i \in \mathcal{I}_{\Sigma}^{\circ}  \},    
\end{equation*}
with
\[
\bfun_0^{(i,i_1)}(0,\xi) =  \left(\phi \circ \pp^{(i_1)}\right)(0, \xi)  \mbox{ and }  \bfun_0^{(i,i_2)}(\xi,0) =  \left(\phi \circ \pp^{(i_2)}\right)(\xi,0),
\]
\[
\bfun_1^{(i,i_1)}(0,\xi) = \frac{\partial_1 \left( \phi \circ \pp^{(i_1)} \right)(0,\xi) + \beta^{(i,i_1)}(\xi) \partial_2 \left( \phi \circ \pp^{(i_1)} \right)(0,\xi) }{\alpha^{(i,i_1)}(\xi)},
\]
and
\[
\bfun_1^{(i,i_2)}(\xi,0) = - \frac{\partial_2 \left( \phi \circ \pp^{(i_2)} \right)(\xi,0) + \beta^{(i,i_2)}(\xi) \partial_1 \left( \phi \circ \pp^{(i_2)} \right)(\xi,0) }{\alpha^{(i,i_2)}(\xi)},
\]
cf.~\cite{FaJuKaTa22}, where each interface curve~$\Sigma^{(i)}$, $i \in \mathcal{I}_{\Sigma}^{\circ}$, is locally parameterised in standard form, see~\ref{appsec:standardform}.
For a function~$\phi \in \mathcal{V}^1$, we denote the equally valued terms $\bfun_0^{(i,i_1)}(0,\xi) = \bfun_0^{(i,i_2)}(\xi,0)$ and $\bfun_1^{(i,i_1)}(0,\xi) = \bfun_1^{(i,i_2)}(\xi,0)$ by the functions $\bfun_0^{(i)}(\xi)$ and $\bfun_{1}^{(i)}(\xi)$, respectively, where $\bfun_0^{(i)}$ represents the trace of the function $\phi$ along the interface curve~$\Sigma^{(i)}$ and $\bfun_1^{(i)}$ describes a specific transversal derivative of the function~$\phi$ across the interface curve~$\Sigma^{(i)}$, cf.~\cite{CoSaTa16}. In case of a boundary curve~$\Sigma^{(i)}  \subset \overline{\Omega^{(i_1)}}$, $i \in \mathcal{I}_{\Sigma}^{\Gamma}$, locally given in standard form, see~\ref{appsec:standardform}, we can equivalently define the functions $\bfun_0^{(i)}$ and $\bfun_{1}^{(i)}$ just as $\bfun_0^{(i)}(\xi)=\bfun_0^{(i,i_1)}(0,\xi)$ and $\bfun_1^{(i)}(\xi)=\bfun_1^{(i,i_1)}(0,\xi)$. 

Due to the dependence of the dimension of $\mathcal{V}^1$ on the initial geometry of the single patch parameterisations~$\pp^{(i)}$, $i \in \mathcal{I}_{\Omega}$, cf.~\cite{KaSaTa17a}, the entire $C^1$-smooth space is, even for simple configurations, hard to study and analyze. Therefore, we consider instead the $C^1$-smooth subspace~$\mathcal{A} \subset \mathcal{V}^1$ introduced in~\cite{FaJuKaTa22}, which is simpler to generate, whose dimension is independent of the initial geometry of the single patch parameterisations~$\pp^{(i)}$, and which still possesses optimal approximation properties as numerically verified in~\cite{FaJuKaTa22}. The $C^1$-smooth subspace~$\mathcal{A}$ is given as   
\[
\mathcal{A} = \{ \phi \in \mathcal{V}^1 \; : \; \bfun_0^{(i)} \in \mathcal{S}_{h}^{p,\reg+1}, \;  \bfun_1^{(i)} \in \mathcal{S}_{h}^{p-1,\reg}, \; i \in \mathcal{I}_{\Sigma}, \mbox{ and } \phi \in C_T^2(\vb{x}^{(i)}),  \; i \in \mathcal{I}_{\chi}  \},
\]
where each edge~$\Sigma^{(i)}$, $i \in \mathcal{I}_{\Sigma}$, is locally parameterised in standard form, see~\ref{appsec:standardform}, and where $\phi \in C_T^2(\vb{x}^{(i)})$ means that the function~$\phi$ is $C^2$-smooth at the vertex~$\vb{x}^{(i)}$ with respect to the tangent plane at $\vb{x}^{(i)}$, cf.~\cite{FaJuKaTa22}.

By requiring $k \geq \frac{4-\reg}{p-\reg-1}$, a possible basis of the space~$\mathcal{A}$ can be given by the set of functions
\[
\Phi = \left( \bigcup_{i \in \mathcal{I}_{\Omega} } \Phi_{\Omega^{(i)}} \right)  \cup \left( \bigcup_{i \in \mathcal{I}_{\Sigma} } \Phi_{\Sigma^{(i)}} \right)  \cup \left( \bigcup_{i \in \mathcal{I}_{\chi} } \Phi_{\vb{x}^{(i)}} \right) 
 \]
with
\[
\Phi_{\Omega^{(i)}} = \{  \phi_{(j_1,j_2)}^{\Omega^{(i)}}  \; : \; j_1,j_2 \in \{2, \ldots, n-3\} \},
\]
\[
\Phi_{\Sigma^{(i)}} = \{ \phi_{(j_1,j_2)}^{\Sigma^{(i)}}  \; : \; j_1 \in \{3-j_2,  \ldots, n_{j_2}-4 +j_2 \} , \; j_2=0,1 \},
\]
and
\[
\Phi_{\vb{x}^{(i)}} = \{ \phi_{(j_1,j_2)}^{\vb{x}^{(i)}}  \; : \; j_1,j_2=0,1,2, \; j_1+j_2 \leq 2 \},
\]
where the sets $\Phi_{\Omega^{(i)}}$, $\Phi_{\Sigma^{(i)}}$ and $\Phi_{\vb{x}^{(i)}}$ collect the basis functions with respect to the individual patches~$\Omega^{(i)}$, edges~$\Sigma^{(i)}$ and vertices~$\vb{x}^{(i)}$, respectively. 
Thereby, all $C^1$-smooth basis functions are locally supported with a support fully contained within the patch~$\Omega^{(i)}$ for a function $\phi_{(j_1,j_2)}^{\Omega^{(i)}}$, $j_1,j_2 \in \{2, \ldots, n-3\}$, with a support in the vicinity of the curve~$\Sigma^{(i)}$ for a function $\phi_{(j_1,j_2)}^{\Sigma^{(i)}}$, $j_1 \in \{3-j_2,  \ldots, n_{j_2}-4 +j_2 \}$, $ j_2=0,1$, and with a support in the vicinity of the vertex~$\vb{x}^{(i)}$ for a function $ \phi_{(j_1,j_2)}^{\vb{x}^{(i)}}$, $j_1,j_2=0,1,2,$ $j_1+j_2 \leq 2$. For the full details of the $C^1$-smooth subspace~$\mathcal{A}$ and of its basis functions, we refer to~\cite{FaJuKaTa22}. A summary of the construction of the basis functions with some implementation details is presented in~\ref{app:construction_space_A}.

\section{Kirchhoff-Love shell analysis of multi-patch structures} \label{sec:numerical_results}
Representing the mid-surface of a multi-patch structured Kirchhoff-Love shell by an AS-$G^1$ multi-patch surface~$\vb{r}$, see Section~\ref{subsec:AS_multi-patch-surface}, and employing the globally $C^1$-smooth multi-patch discretisation space~\cite{FaJuKaTa22}, recalled in Section~\ref{subsec:specific_space}, we can directly apply the Kirchhoff-Love shell formulation~\cite{kiendl-bletzinger-linhard-09}, introduced in Section~\ref{sec:KLShell}, to the individual surfaces patches~$\vb{r}^{(i)}$ of the multi-patch mid-surface~$\vb{r}$. This leads to a novel, simple isogeometric method for the analysis of Kirchhoff-Love shells composed of multiple patches with possibly extraordinary vertices. The performance of this method will be tested in this section on the basis of a series of benchmark problems. In \cref{subsec:hyperboloid} the linear Kirchhoff-Love shell equation will be solved for a hyperboloid shell, adding a hole for the benchmark in \cref{subsec:hyperboloid_hole}. In \cref{subsec:L-shape}, post-buckling analysis of a clamped corner piece will be performed. In \cref{subsec:L-shape_hole}, holes are added to this geometry. For the geometric and material properties of the considered shells we will denote by $L, W$ and $a$ the lengths, by $t$ the thickness, by $E$ the Young's modulus, and by $\nu$ the Poisson's ratio. In all examples, we will perform uniform h-refinement for the convergence analysis. In cases where the penalty method from \cite{Herrema2019} is employed, a penalty parameter of $\alpha=10^3$ will be used, unless another value will be required to produce good convergence rates and stress fields.

\subsection{Hyperboloid shell}\label{subsec:hyperboloid}
\begin{figure}
\centering
\begin{subfigure}[t]{\linewidth}
    \centering
    \begin{minipage}{0.4\linewidth}
    \resizebox{\linewidth}{!}{\includegraphics{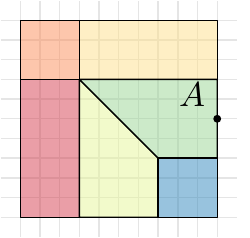}}
    \end{minipage}
    \hfill
    \begin{minipage}{0.55\linewidth}
    \centering
    \includegraphics[width=\linewidth]{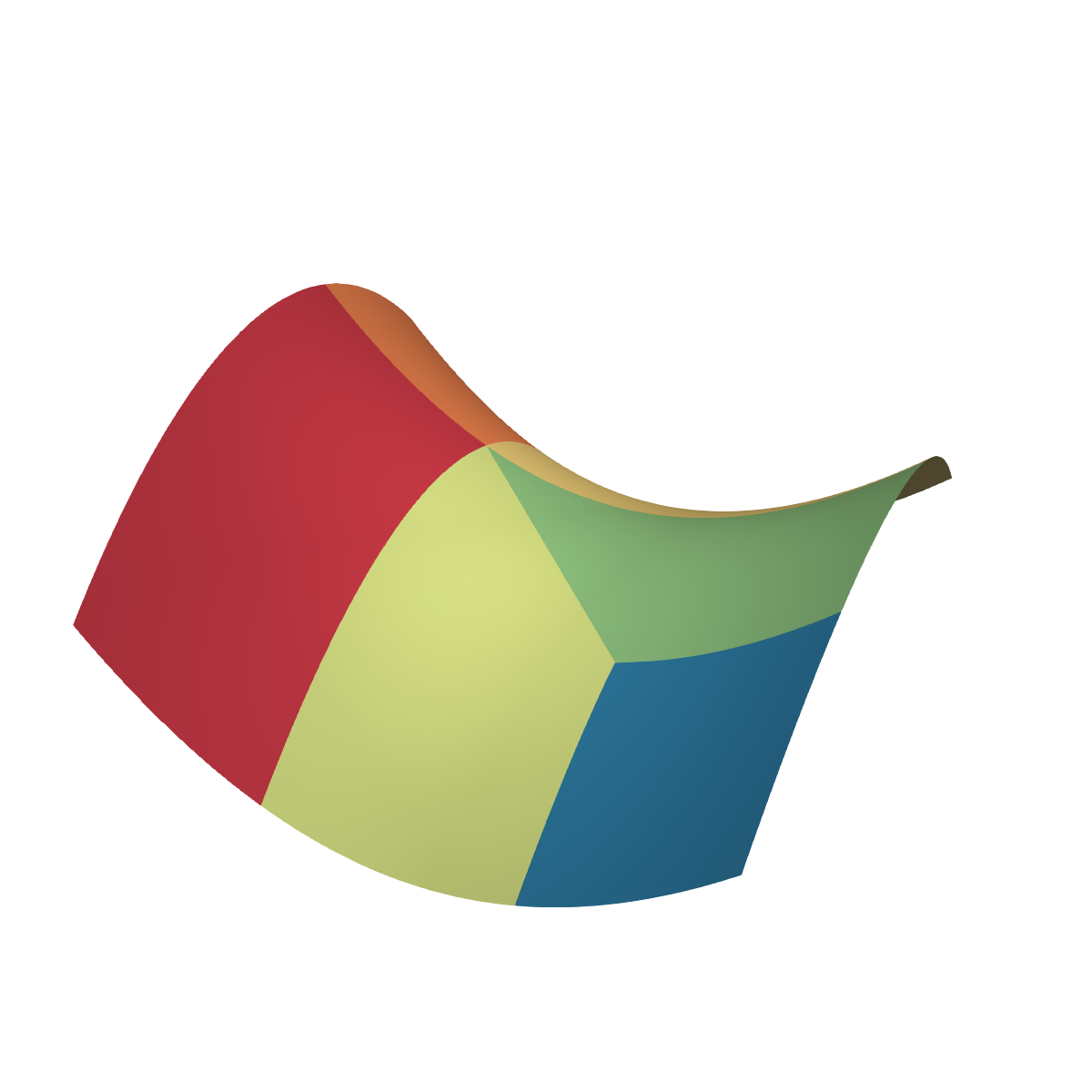}
    \end{minipage}
    \caption{Geometry 1}
    \label{fig:hyperboloid1}
\end{subfigure}

\begin{subfigure}[t]{\linewidth}
    \centering
    \begin{minipage}{0.4\linewidth}
    \resizebox{\linewidth}{!}{
    \includegraphics{Figures/6p_planar.pdf}
    }
    \end{minipage}
    \hfill
    \begin{minipage}{0.55\linewidth}
    \centering
    \includegraphics[width=\linewidth]{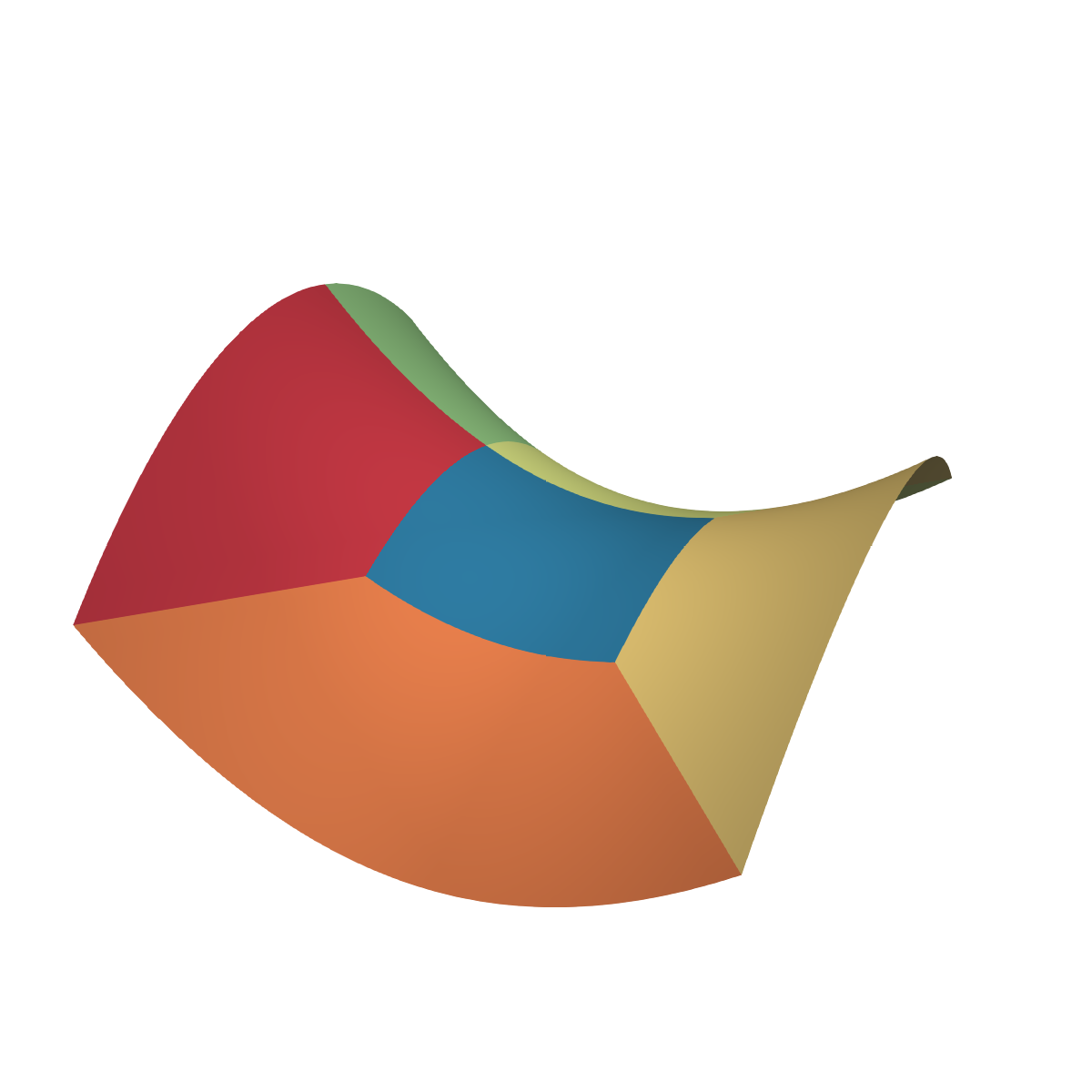}
    \end{minipage}
    \caption{Geometry 2}
    \label{fig:hyperboloid2}
\end{subfigure}
\caption{AS-$G^1$ mid-surfaces of hyperboloid multi-patch shell geometries consisting of 6 patches in two different configurations. For both subfigures (a) and (b) the left figure represents the parametric domain and the right figure gives and impression of the resulting geometry. The reference point~$A$ is used in the numerical tests, cf. \cref{fig:hyperboloid_results}.}
\label{fig:hyperboloid}
\end{figure}

We represent the mid-surface of a hyperboloid shell,  given by the hyperbolic surface
\begin{equation} \label{eq:surface1}
 \surf(\xi^1,\xi^2) = \begin{bmatrix}
                   \xi^1 & \xi^2 & (\xi^1)^2-(\xi^2)^2
                  \end{bmatrix}^T
\end{equation}
for the parameter domain~$D=[-L/2, L/2]^2$, by the two different 6-patch AS-$G^1$ multi-patch surfaces ~$\pp_{\ell}$, $\ell=1,2$, visualized in \cref{fig:hyperboloid}. The multi-patch parameterisation from \cref{fig:hyperboloid1} is selected since it has two interior extraordinary vertices, i.e. vertices where 3 or more than 4 patches join. The geometry from \cref{fig:hyperboloid2} is selected due to its interior extraordinary vertices and due to the fact that it has boundary vertices, where two patches join $C^0$ in the south-east, south west and north east corners of \cref{fig:hyperboloid2}; an aspect where the method from \cite{Toshniwal2017} would smoothen the geometry.

For the design of the multi-patch surfaces in \cref{fig:hyperboloid} we follow the strategy from \cite[Example~3]{FaJuKaTa22}. Let $\widetilde{\pp}_{\ell}$, $\ell=1,2$, be the bilinearly parameterised multi-patch parameterisations shown in \cref{fig:hyperboloid1,fig:hyperboloid2} (left), respectively, which represent in each case the parameter domain~$D$. We then generate the multi-patch surfaces~$\pp_{\ell}$, $\ell=1,2$, just by choosing the single surface patch parameterisations $\pp_{\ell}^{(i)}$ as~$\pp_{\ell}^{(i)} = \surf \circ \widetilde{\pp}_{\ell}^{(i)}$, where $\widetilde{\pp}_{\ell}^{(i)}$ are the individual patch parameterisations of the parameterisations~$\widetilde{\pp}_{\ell}$, which implies that $\pp_{\ell}^{(i)} \in (\mathcal{S}_h^{\vb{p}, \breg})^3$  with $\vb{p}=(2,2)$, $\breg=(\infty, \infty)$ and $h=1$. The two resulting multi-patch surfaces~$\pp_{\ell}$, $\ell=1,2$, are AS-$G^1$ by construction, since the surface~$\surf$ is just a bivariate polynomial surface patch, and the bilinear parameterisations $\widetilde{\pp}_{\ell}$, $\ell=1,2$, are trivially AS-$G^1$.

\begin{figure}
    \centering
    \includegraphics[width=0.87\linewidth]{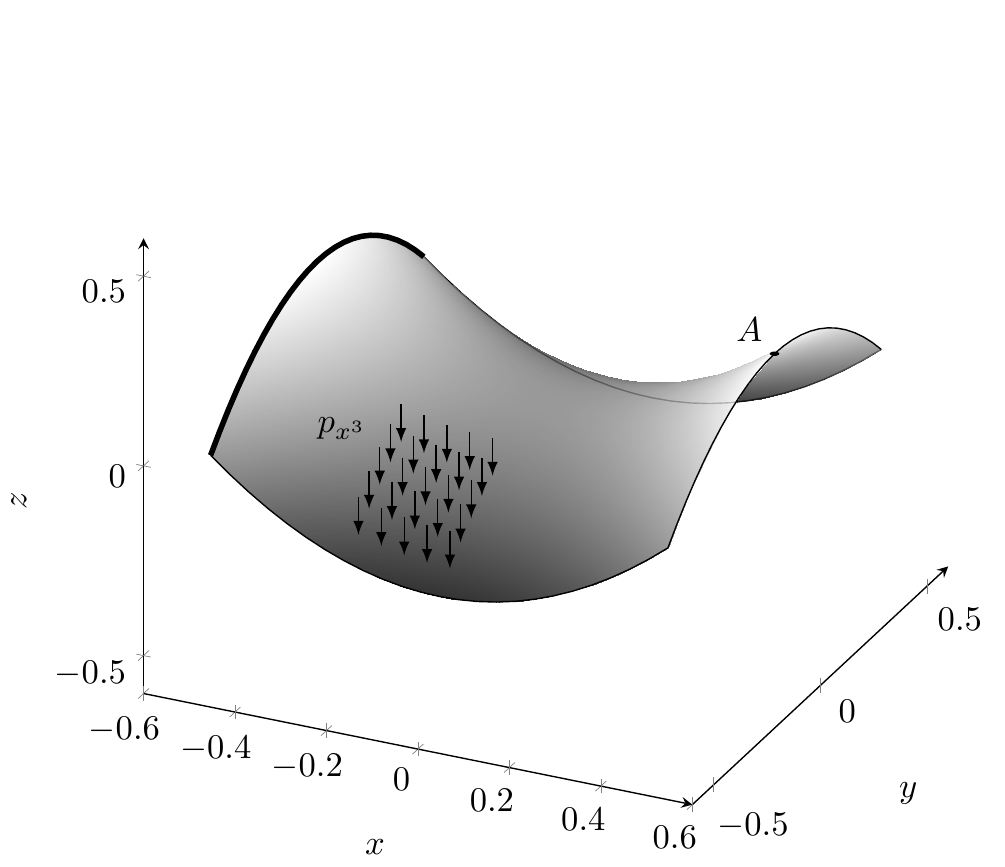}
    \caption{Set-up of the clamped hyperboloid problem from \cite{Belytschko1985}. The displacements and rotations are fixed on the left boundary represented by the thick line. The distributed load $p_{x^3}$ is applied over the whole domain in vertical direction. Furthermore, the point $A$ is used for reference.}
    \label{fig:hyperboloid_setup}
\end{figure}

The benchmark problem is based on \cite{Belytschko1985}, and an overview is depicted in \cref{fig:hyperboloid_setup}. The geometric and material properties of the hyperbolic multi-patch shells $\vb{y}_{\ell}$, $\ell=1,2$, are determined by the length $L=1\:[\text{m}]$, thickness $t=0.01\:[\text{m}]$, Young's modulus $E=2 \cdot 10^{11}\:[\text{N}/\text{m}^2]$ and the Poisson's ratio $\nu = 0.3\:[-]$. The shell is fully clamped on the west side at $x^1 = -L/2$, which means that all displacements and rotations are fixed weakly \cite{Herrema2019}. We apply a uniformly distributed vertical load with magnitude $p_{x^3} = -8000 \cdot t\:[\text{N}/\text{m}^2]$ and perform linear analysis. We study the vertical displacement at point $A=(x^1,x^2,x^3) = (L/2,0,L^4/4)$ as in \cite{Belytschko1985} as well as the strain energy norm $B=\frac{1}{2}\vb{u}^T K \vb{u}$, where $\vb{u}$ is the solution displacement vector and $K$ is the linear stiffness matrix. Numerical tests are performed on the two multi-patch shells~$\vb{y}_{\ell}$, $\ell=1,2$, possessing the mid-surfaces~$\vb{r}_{\ell}$, for polynomial degrees $p = 3,4, 5$, and maximal possible regularity $\reg=p-2$ for the $C^1$ spline space~$\mathcal{A}$. These geometries are presented in \cref{fig:hyperboloid}.

\begin{figure}
\centering
    \begin{tabular}{cc}
     \includegraphics[width=0.44\linewidth]{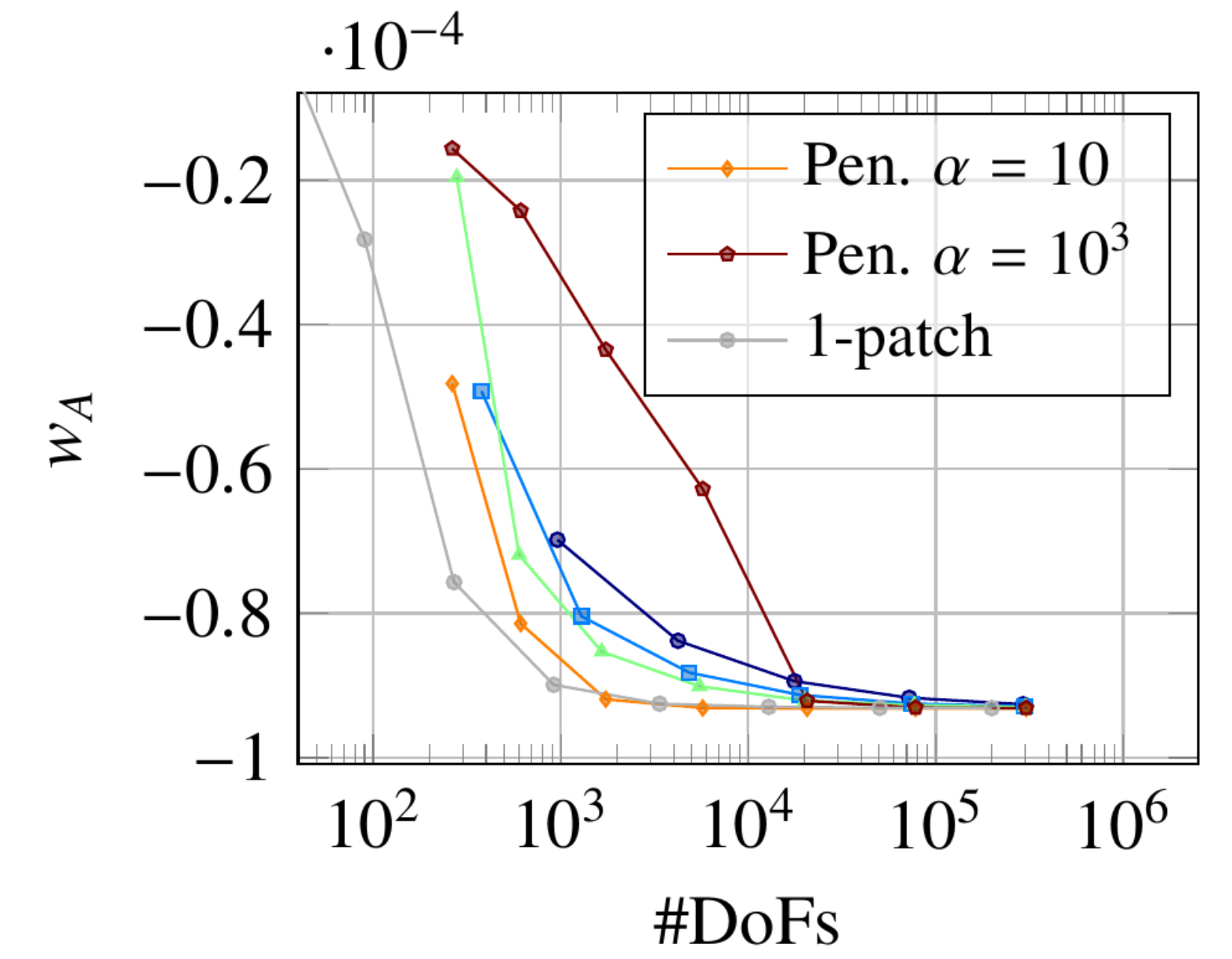} & 
     \includegraphics[width=0.44\linewidth]{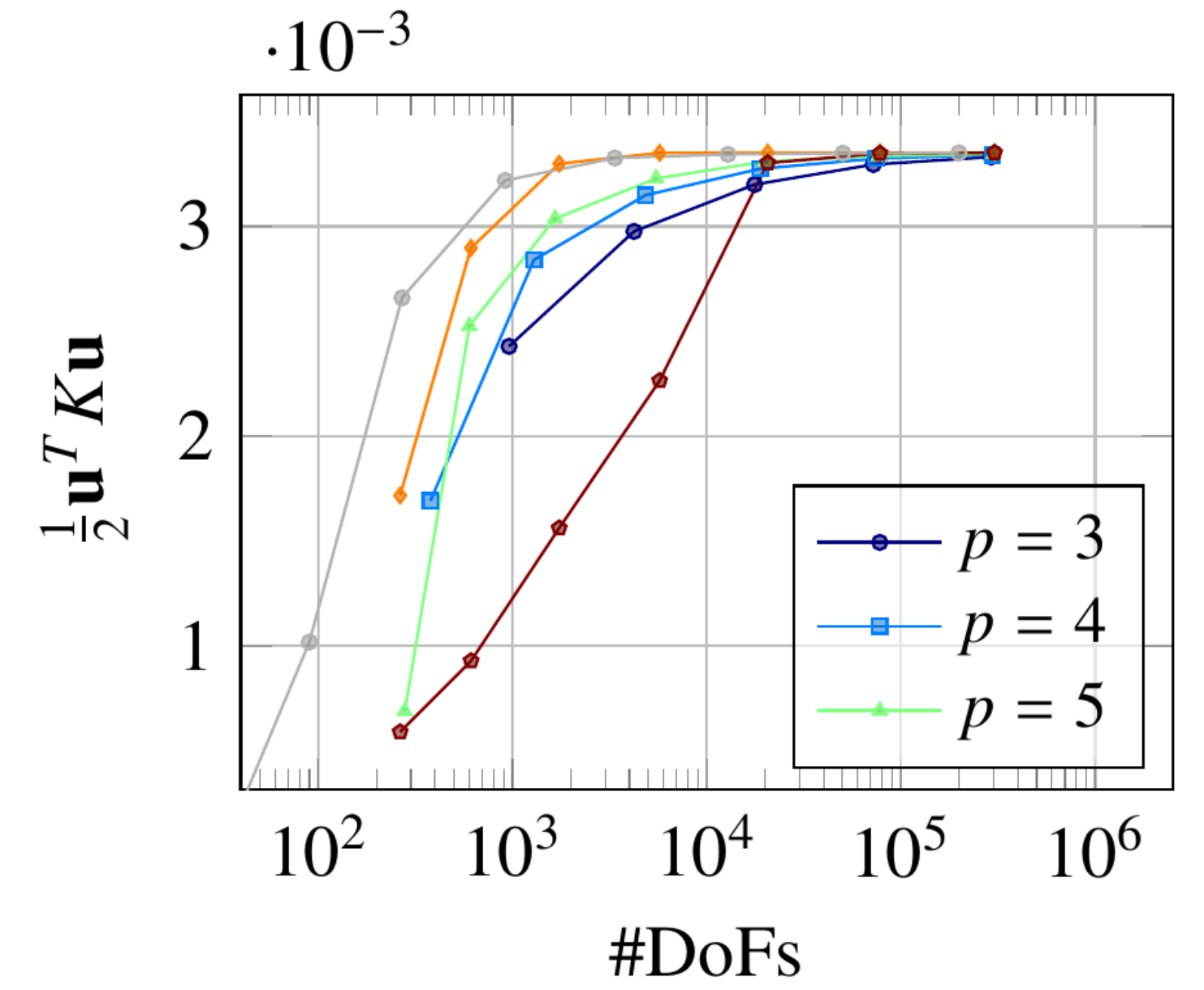} \\
     \multicolumn{2}{c}{(a) Geometry 1} \\
     \includegraphics[width=0.44\linewidth]{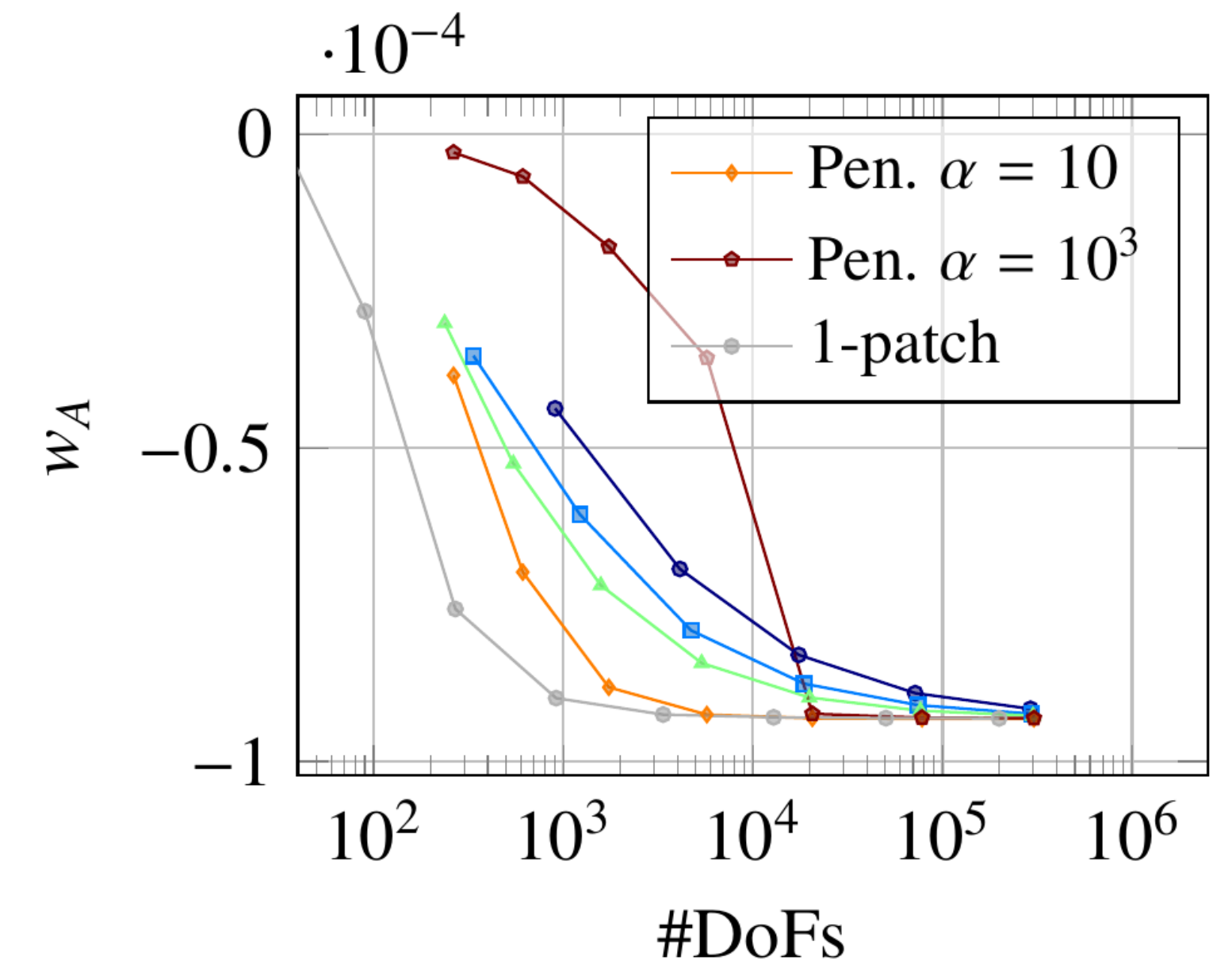} & 
     \includegraphics[width=0.44\linewidth]{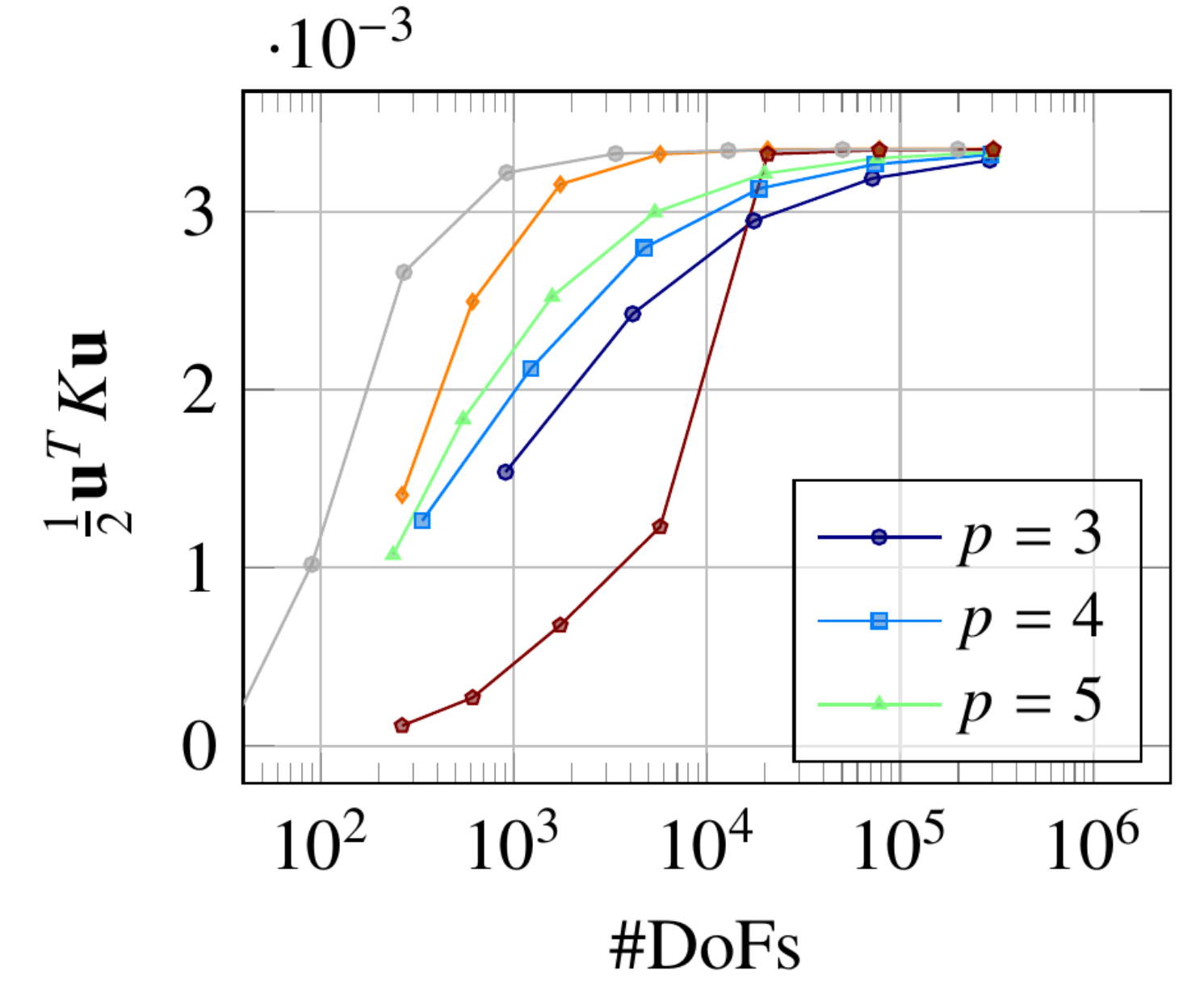} \\
     \multicolumn{2}{c}{(b) Geometry 2} 
    \end{tabular}
    \caption{Results for the 6-patch hyperboloid geometries from \cref{fig:hyperboloid} subject to a uniform vertical load. The left plots represent the vertical displacement in the reference point $A$ ($w_A$) and the right plots represent the strain energy norm ($\vb{u}^T K \vb{u}$), both against the number of degrees of freedom (\#DoFs). For the AS-$G^1$ constructions, the degree $p$ is varied between $3$ and $5$ with maximum regularity $\reg=p-2$. The results for the single-patch case as well as for the penalty method are obtained for degree $p=3$ and regularity $\reg=1$ and in case of the penalty method using the penalty parameters $\alpha=10^1$ and $\alpha=10^3$.}
    \label{fig:hyperboloid_results}
\end{figure}

The results for the fully clamped hyperboloid shells with the two different AS-$G^1$ multi-patch configurations for the mid-surfaces are presented in \cref{fig:hyperboloid_results}. The convergence plots include for comparison the results for a single patch geometry and for a penalty coupling with penalty parameters $\alpha=10^1, 10^3$, see \cite{Herrema2019}, with degree $p=3$ and regularity $\reg=1$. The penalty parameter of $\alpha=10^1$ is selected for comparison purposes, given its convergence to the single patch solution. It can be seen that the AS-$G^1$ multi-patch results show in the point-wise displacement as well as in terms of the energy norm convergence towards the single-patch results, but with a lower rate as the results for the single-patch or for the penalty coupling case.

\begin{figure}
\centering
\begin{subfigure}[t]{0.9\linewidth}
    \centering
        Von Mises stress [MPa]
        \includegraphics[width=0.8\linewidth,trim=0 0 0 1105,clip]{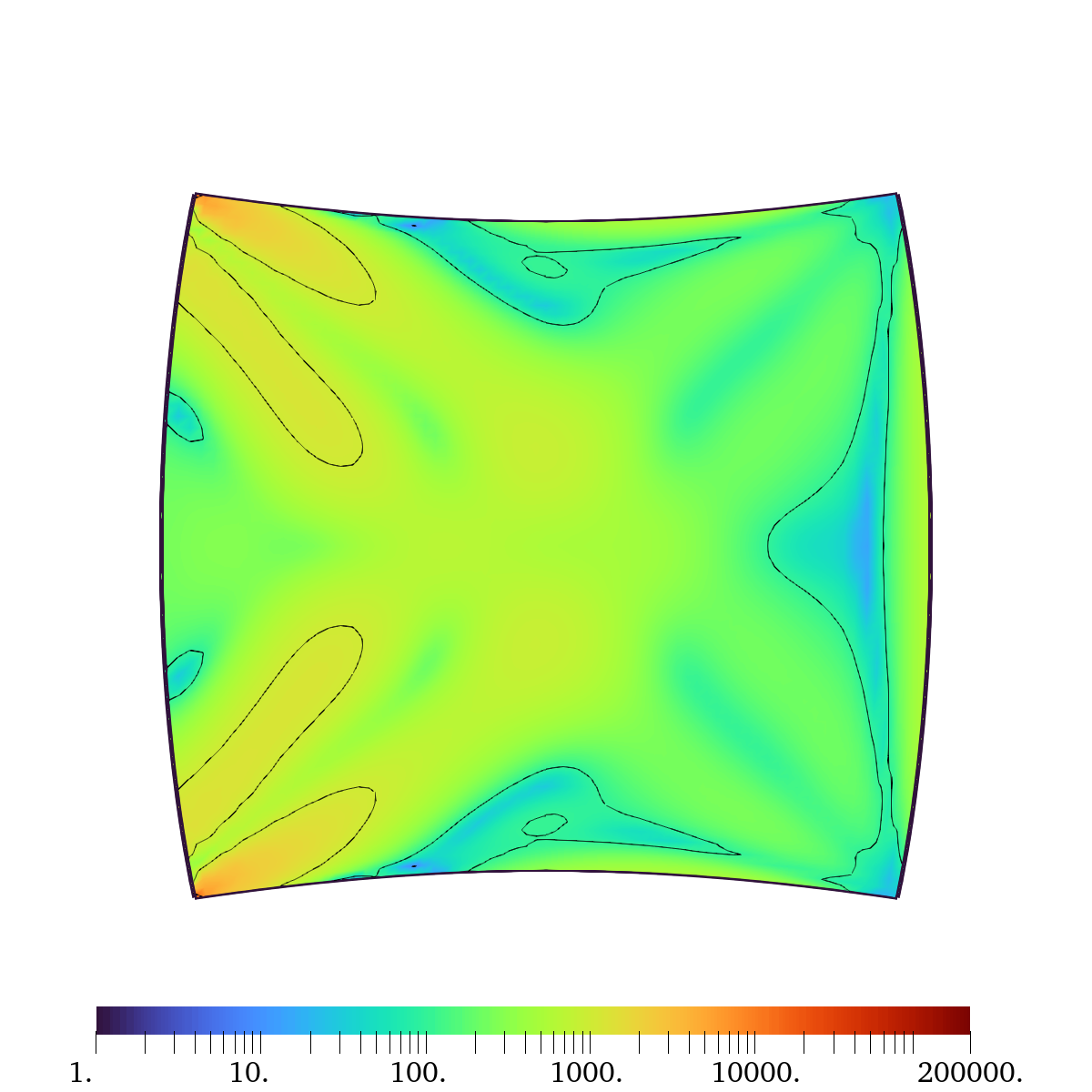}
\end{subfigure}

\begin{subfigure}[t]{0.45\linewidth}
    \centering
        \includegraphics[width=0.8\linewidth]{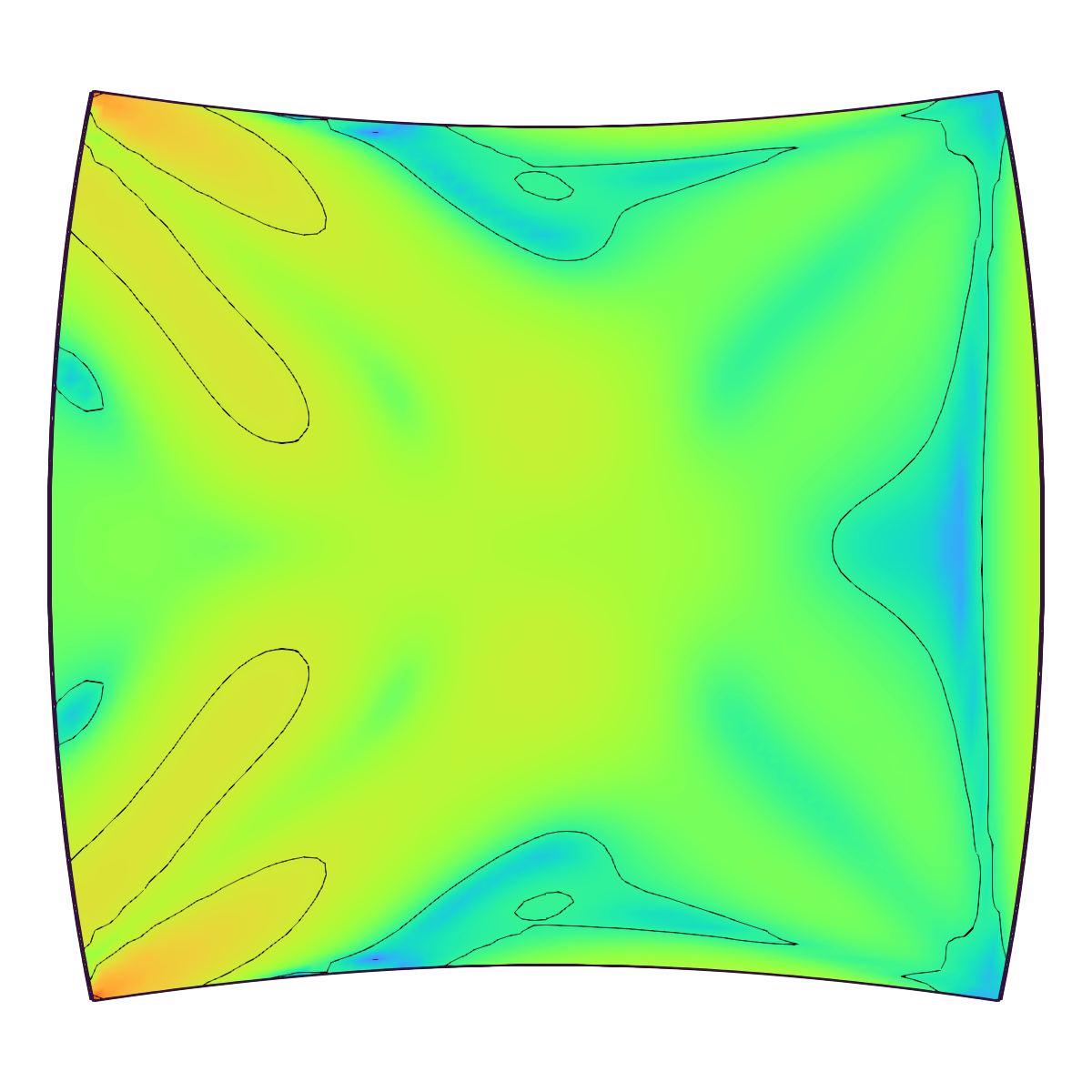}
    \caption{Single-patch geometry.}
\end{subfigure}
\hfill
\begin{subfigure}[t]{0.45\linewidth}
\end{subfigure}

\begin{subfigure}[t]{0.45\linewidth}
    \centering
        \includegraphics[width=0.8\linewidth]{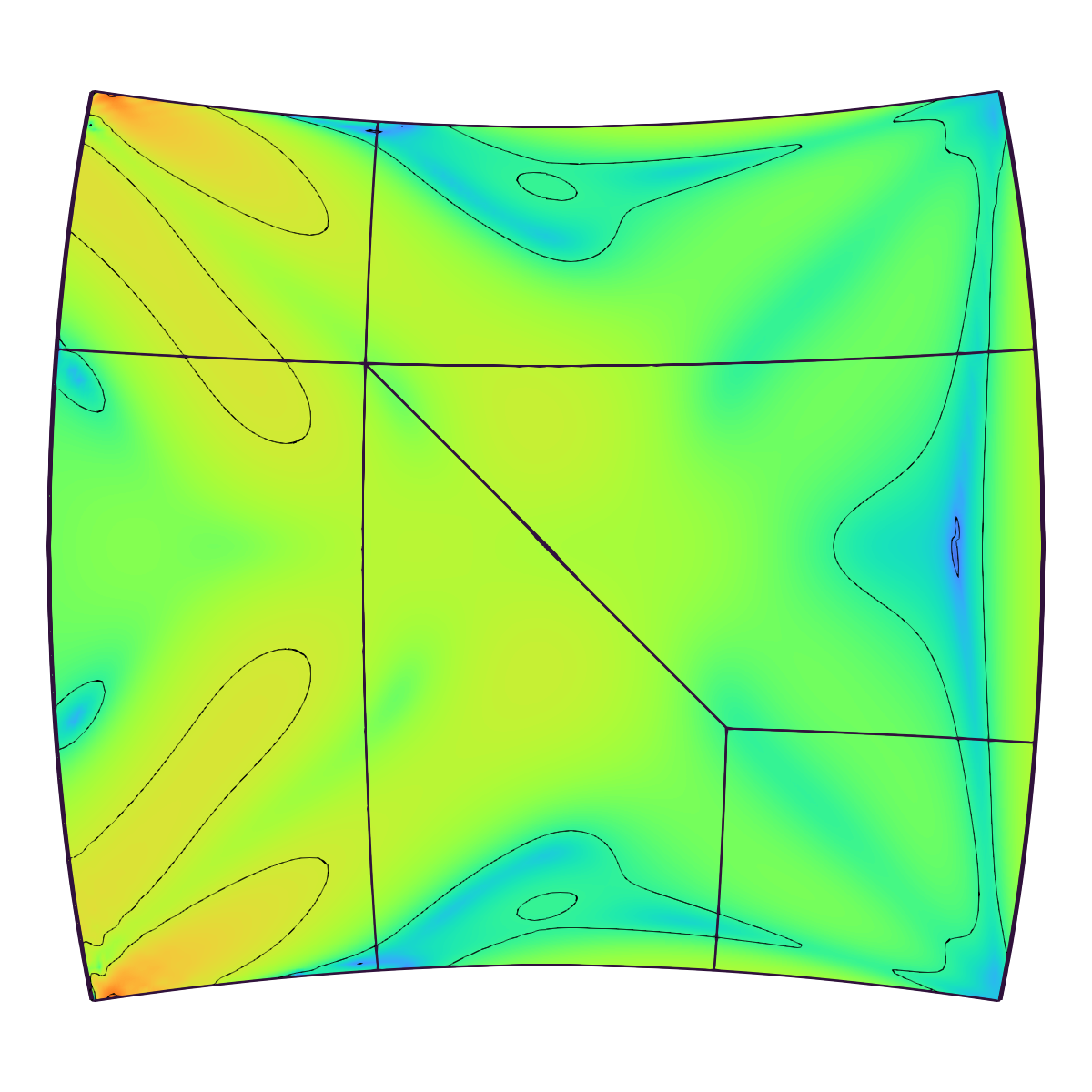}
    \caption{Geometry 1, AS-$G^1$ method.}
\end{subfigure}
\hfill
\begin{subfigure}[t]{0.45\linewidth}
    \centering
        \includegraphics[width=0.8\linewidth]{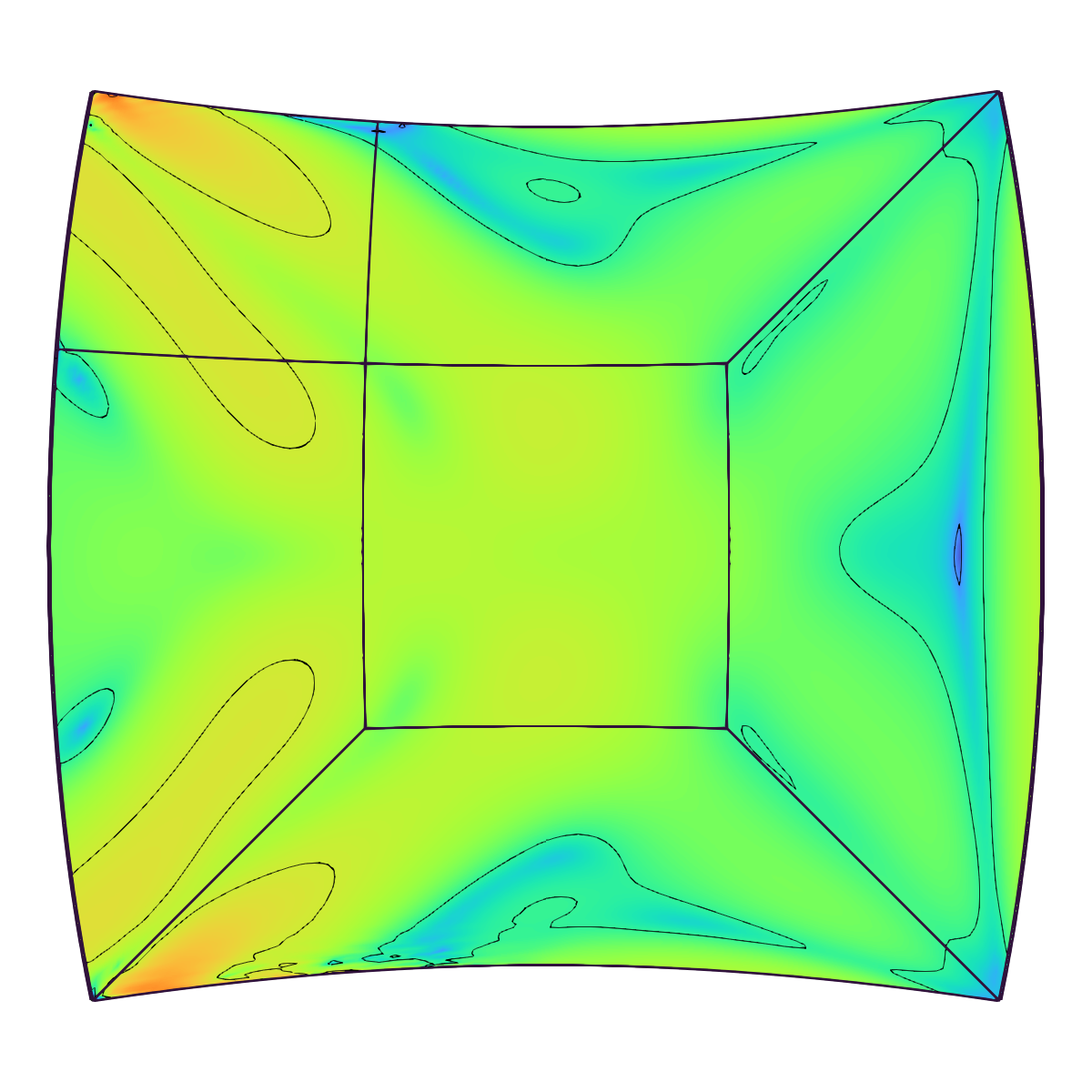}
    \caption{Geometry 2, AS-$G^1$ method.}
\end{subfigure}

\begin{subfigure}[t]{0.45\linewidth}
    \centering
        \includegraphics[width=0.8\linewidth]{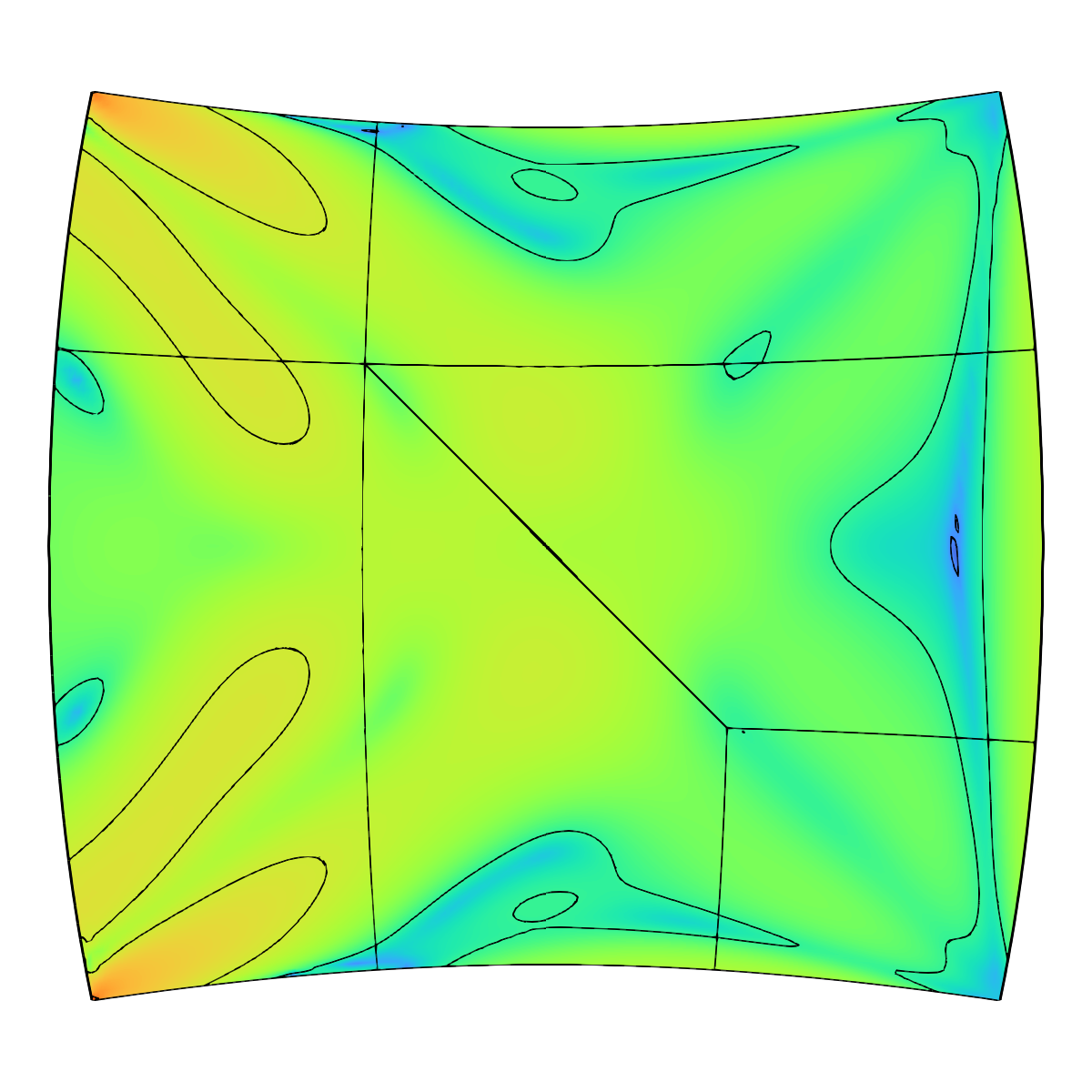}
    \caption{Geometry 1, penalty method.}
\end{subfigure}
\hfill
\begin{subfigure}[t]{0.45\linewidth}
    \centering
        \includegraphics[width=0.8\linewidth]{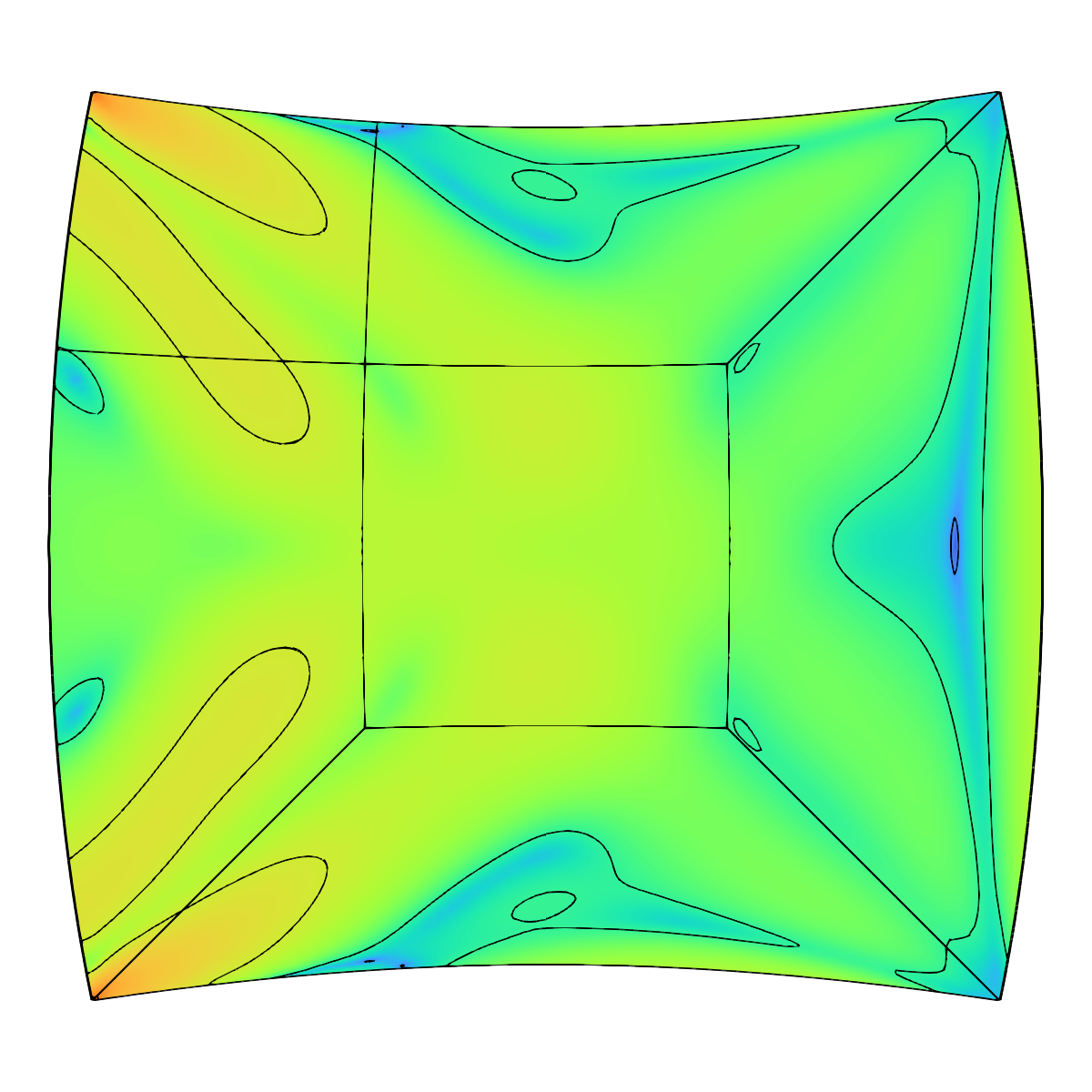}
    \caption{Geometry 2, penalty method.}
\end{subfigure}
\caption{Von Mises membrane stress fields for the single-patch and multi-patch hyperboloid geometries in \cref{fig:hyperboloid} on a $16\times16$ element mesh per patch for degree $p=4$ and regularity $\reg=2$. The penalty parameter for the penalty method results is $\alpha=10^1$. All stress fields are plotted on the same colour scale. The contour lines are provided for values $10$, $10^2$, $10^3$, $10^4$ and $10^5\:[\text{MPa}]$.} 
\label{fig:hyperboloid_stress}
\end{figure}

Furthermore, the Von Mises membrane stress is computed for the single-patch hyperboloid as well as for the considered multi-patch geometries. The results are obtained for degree $p=4$ and regularity $\reg=2$ with $16\times16$ elements per patch. Since the regularity is equal to 2, it is expected that the stresses are $C^1$ continuous in the interior and $C^0$ continuous over the patch boundaries, as they are derivatives of the $C^2$ interior basis which is $C^1$ over the patch boundaries by construction. The results obtained by the penalty method are also given for degree $p=4$, regularity $\reg=2$ and for the penalty parameter $\alpha=10^1$. \cref{fig:hyperboloid_stress} verifies that the stress fields are smooth and well represented by our AS-$G^1$ multi-patch construction compared to the single-patch stress field and the field obtained by the penalty method. For the AS-$G^1$ geometry, small artifacts can be observed on the top and bottom boundaries for the stress field. This could be simply solved in future by using for each boundary edge e.g. the entire space instead of the subspace which has been selected in~\cite{FaJuKaTa22} to have a simple and uniform construction. Furthermore, it should be noted that the results of the penalty method depend on the choice of a coupling parameter, as shown in \cref{fig:hyperboloid_results}, whereas the present method is parameter-free.

\subsection{Hyperboloid shell with a hole}\label{subsec:hyperboloid_hole}

\begin{figure}
    \centering
    \includegraphics[width=0.87\linewidth]{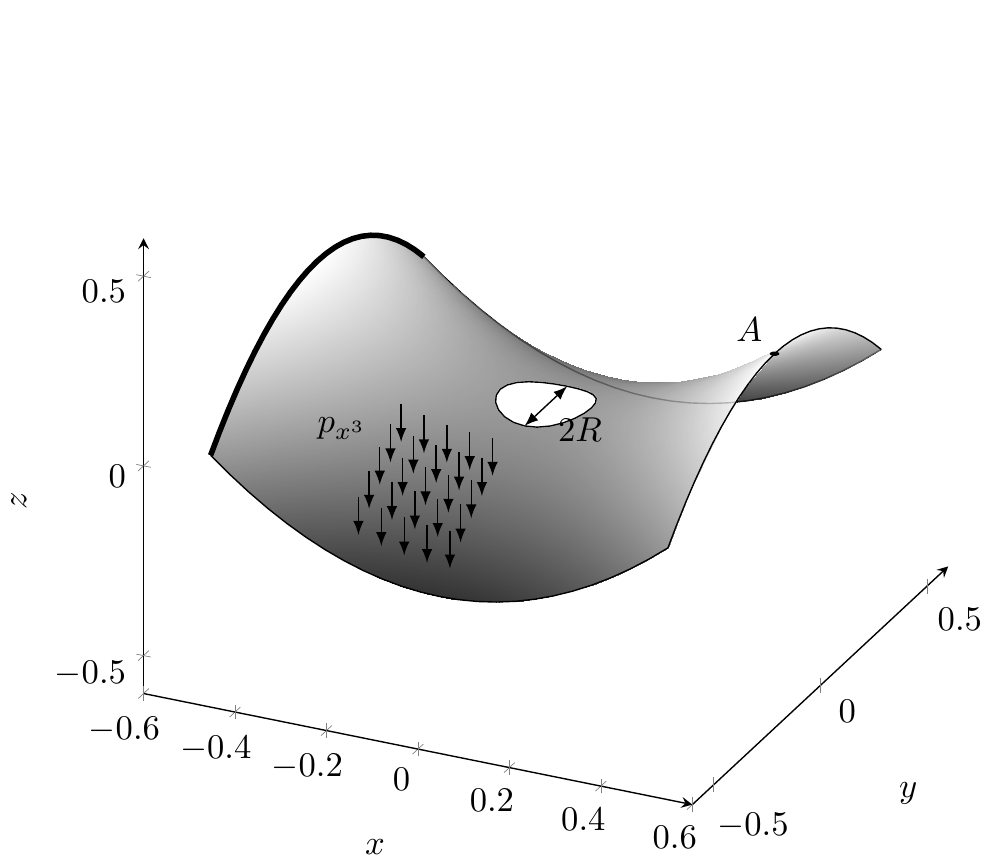}
    \caption{Set-up of the clamped hyperboloid with a hole of radius $R=0.15$. The displacements and rotations are fixed on the left boundary represented by the thick line. The distributed load $p_{x^3}$ is applied over the whole domain in vertical direction. Furthermore, the point $A$ is used for reference.}
    \label{fig:hyperboloid_hole_setup}
\end{figure}

\begin{figure}[htbp]
\centering
\begin{subfigure}[t]{0.45\linewidth}
    \centering
    \includegraphics[width=\linewidth]{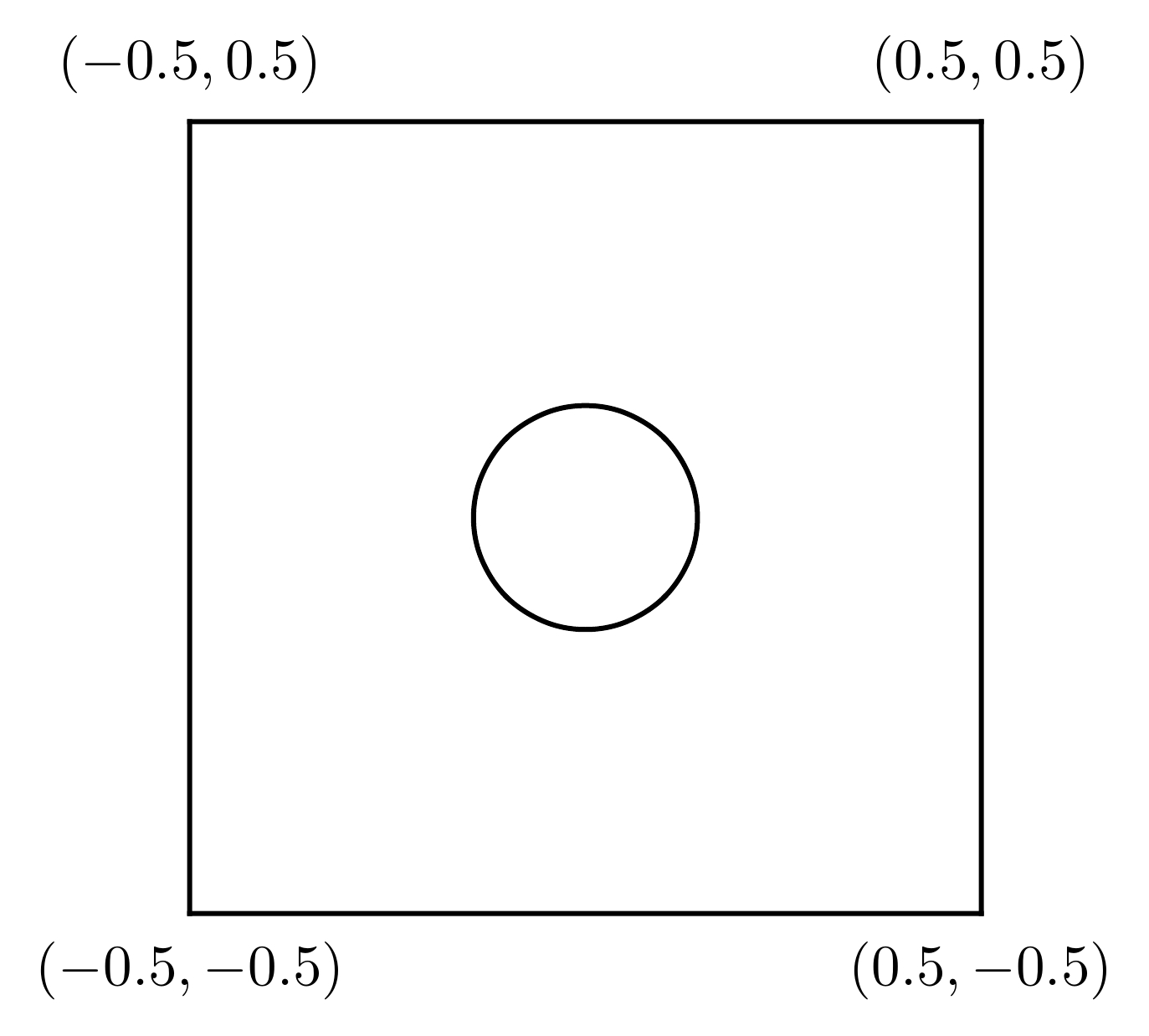}
    \caption{Trimmed parameter domain}
    \label{fig:hyperboloid_trimmed}
\end{subfigure}
\hfill
\begin{subfigure}[t]{0.45\linewidth}
    \centering
    \includegraphics[width=\linewidth]{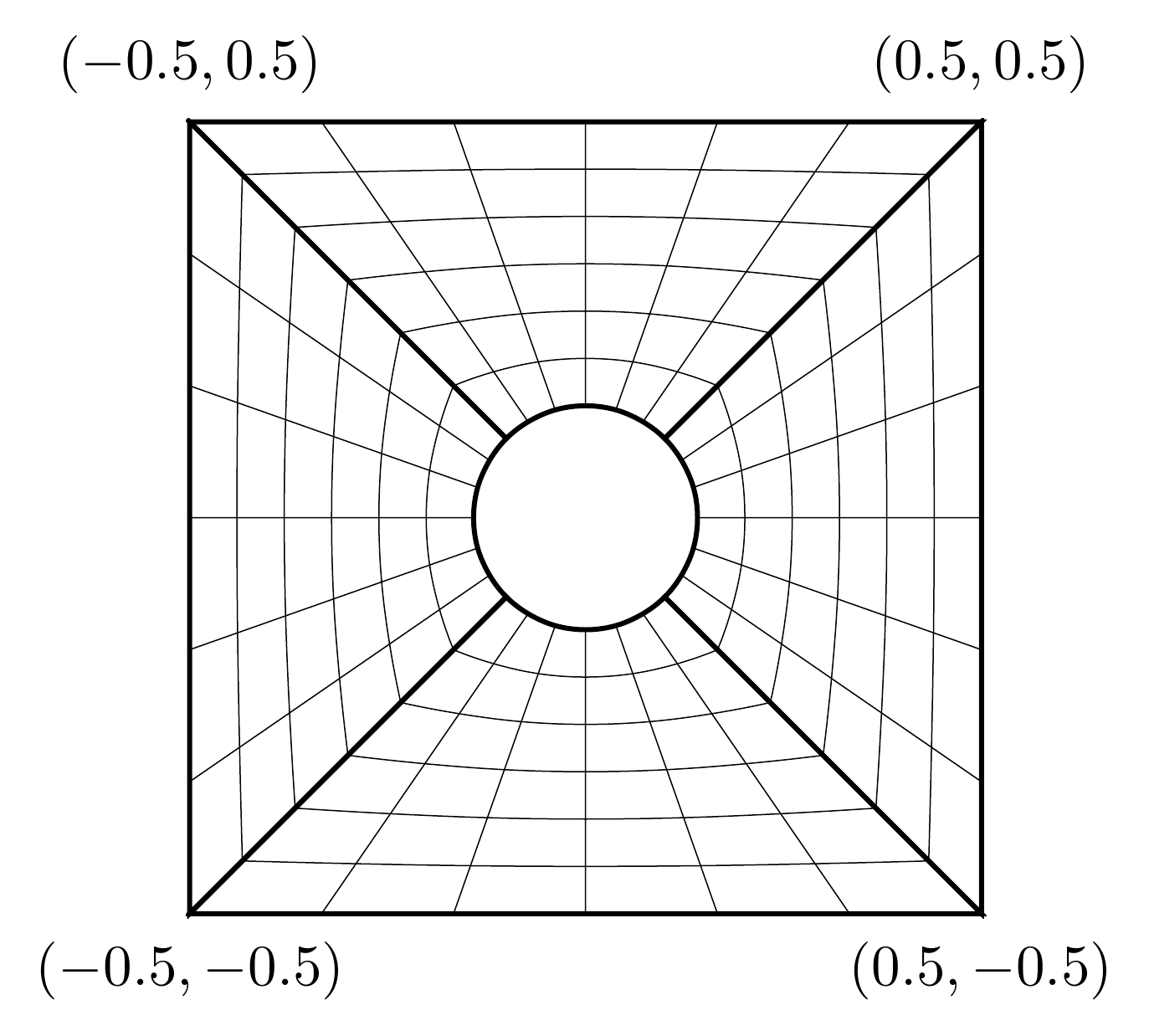}
    \caption{Untrimmed parameter domain}
    \label{fig:hyperboloid_untrimmed}
\end{subfigure}

\begin{subfigure}[t]{0.55\linewidth}
    \centering
    \includegraphics[width=\linewidth]{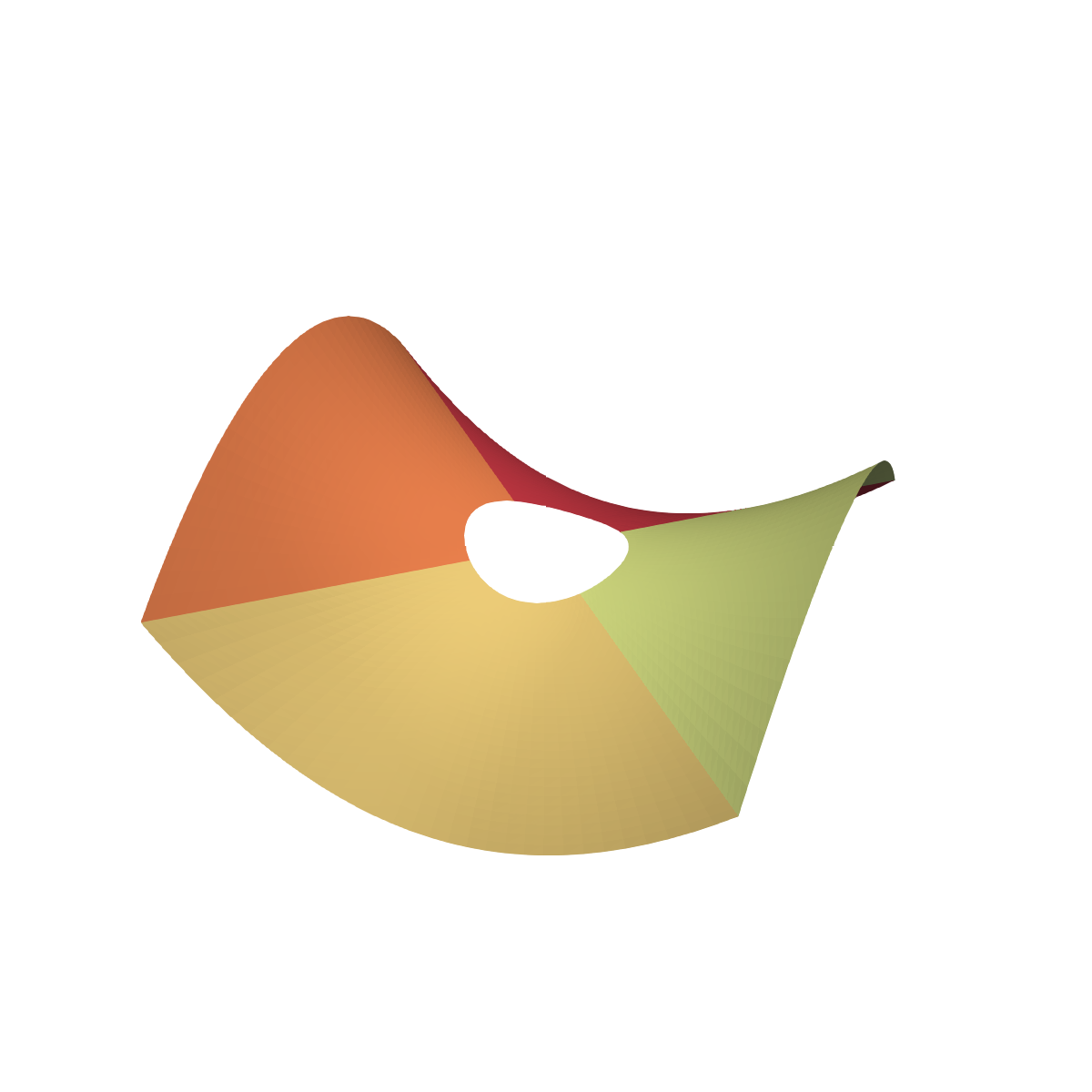}
    \caption{AS-$G^1$ 4-patch geometry}
\end{subfigure}
\caption{Multi-patch hyperboloid shell with a hole. (a): Trimmed parameter domain. (b): Untrimmed parameter domain given by an AS-$G^1$ 4-patch parameterisation. (c) AS-$G^1$ 4-patch surface of the hyperboloid shell with hole with the reference point $A$ for the numerical tests, cf. \cref{fig:hyperboloid_hole_results}.
}
\label{fig:hyperboloid_hole}
\end{figure}

We consider as in the previous example in \cref{subsec:hyperboloid} a hyperboloid shell but now with a hole, see \cref{fig:hyperboloid_hole_setup}. Although this geometry does not possess extraordinary vertices, it cannot be represented by a regular single-patch surface due to the presence of vertices on the outer boundary, where two patches join $C^0$, and of the presence of the smooth inner boundary at the same time. The AS-$G^1$ 4-patch geometry representing the mid-surface of the shell is constructed via trimming from the hyperbolic surface~$\widehat{\vb{r}}$, given in~\eqref{eq:surface1}, by following the strategy from~\cite[Example~2]{FaJuKaTa22}. For this, we first trim the parameter domain $D=[-L/2, L/2]^2$ of the hyperbolic surface by cutting out an approximated disk with radius $R=0.15$ whose boundary is a B-spline curve of degree~$2$, see \cref{fig:hyperboloid_trimmed}. The associated untrimmed parameter domain (see \cref{fig:hyperboloid_untrimmed}) is then described via a planar AS-$G^1$ 4-patch parameterisation~$\widetilde{\pp}$, consisting of individual patch parameterisations $\widetilde{\pp}^{(i)}$, $i =1,\ldots, 4$, with $\widetilde{\pp}^{(i)} \in \left(\mathcal{S}_{h}^{\vb{p}, \breg}\right)^2$, $\vb{p}=(2,2)$, $\breg=(1,1)$ and $h=\frac{1}{4}$, see \cref{fig:hyperboloid_hole} (b). Via $\pp^{(i)} = \surf \circ \widetilde{\pp}^{(i)}$, $i =1,\ldots, 4$, we obtain surface parameterisations $\pp^{(i)} \in \left(\mathcal{S}_h^{\vb{p}, \breg}\right)^3$, with $\vb{p}=(4,4)$, $\breg=(1,1)$ and $h=\frac{1}{4}$, which finally form together the AS-$G^1$ 4-patch surface~$\pp$, shown in \cref{fig:hyperboloid_untrimmed}, representing the mid-surface of the desired hyperboloid shell with hole. Note that this geometry cannot be used by the method from \cite{Toshniwal2017}, since it contains vertices on the boundary, where two patches join~$C^0$.  

\begin{figure}
    \centering
    \begin{tabular}{cc}
     \includegraphics[width=0.44\linewidth]{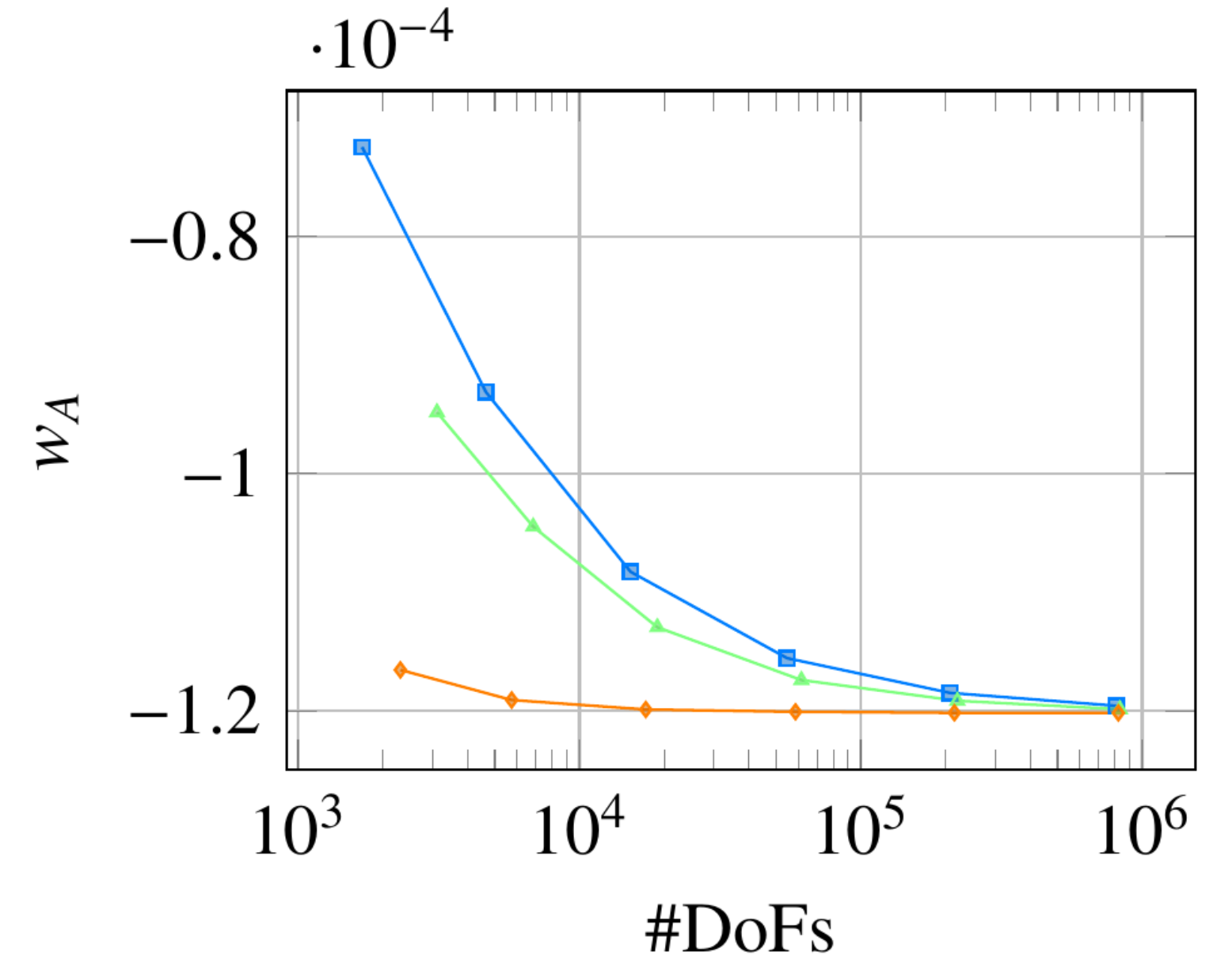} &
     \includegraphics[width=0.44\linewidth]{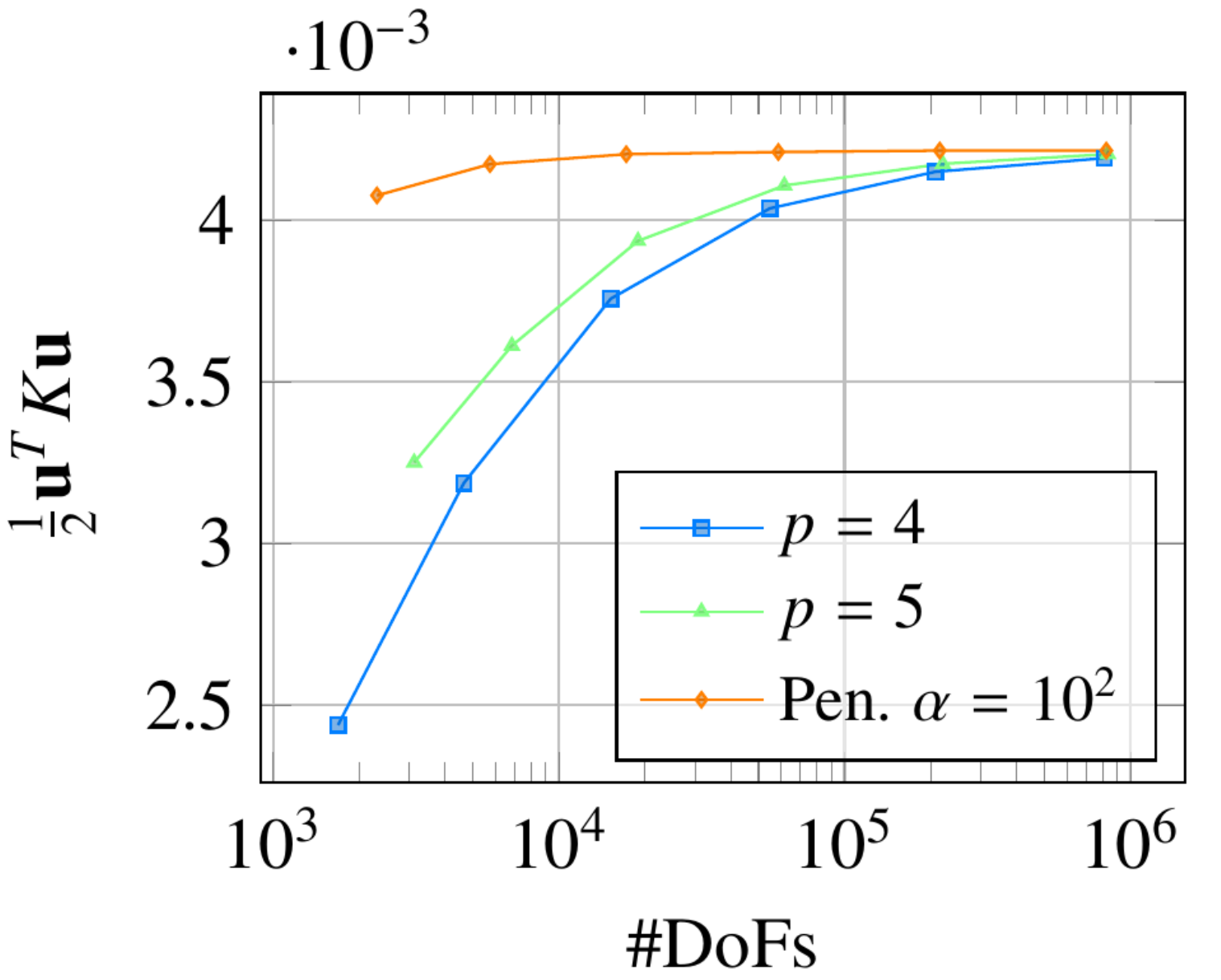} 
    \end{tabular}
    \caption{Results for the 4-patch hyperboloid with a hole from \cref{fig:hyperboloid_hole} subject to a uniform vertical load. The left plot represents the vertical displacement in the reference point $A$ ($w_A$) and the right plot represents the strain energy norm ($\vb{u}^T K \vb{u}$), both against the number of degrees of freedom (\#DoFs). For the AS-$G^1$ constructions, the degree $p$ is varied between $4$ and $5$ with maximum regularity $\reg=p-2$. The results for the penalty method is obtained for degree $p=4$ and regularity $\reg=2$ using the penalty parameter $\alpha=10^2$.}
    \label{fig:hyperboloid_hole_results}
\end{figure}

\begin{figure}
    \centering
    \begin{subfigure}[t]{0.9\linewidth}
        \centering
            Von Mises stress [MPa]
            \includegraphics[width=0.8\linewidth,trim=0 0 0 1105,clip]{Figures/1p_hyperboloid_stress_colorbar.png}
    \end{subfigure}

    \begin{subfigure}[t]{0.45\linewidth}
        \centering
            \includegraphics[width=0.8\linewidth]{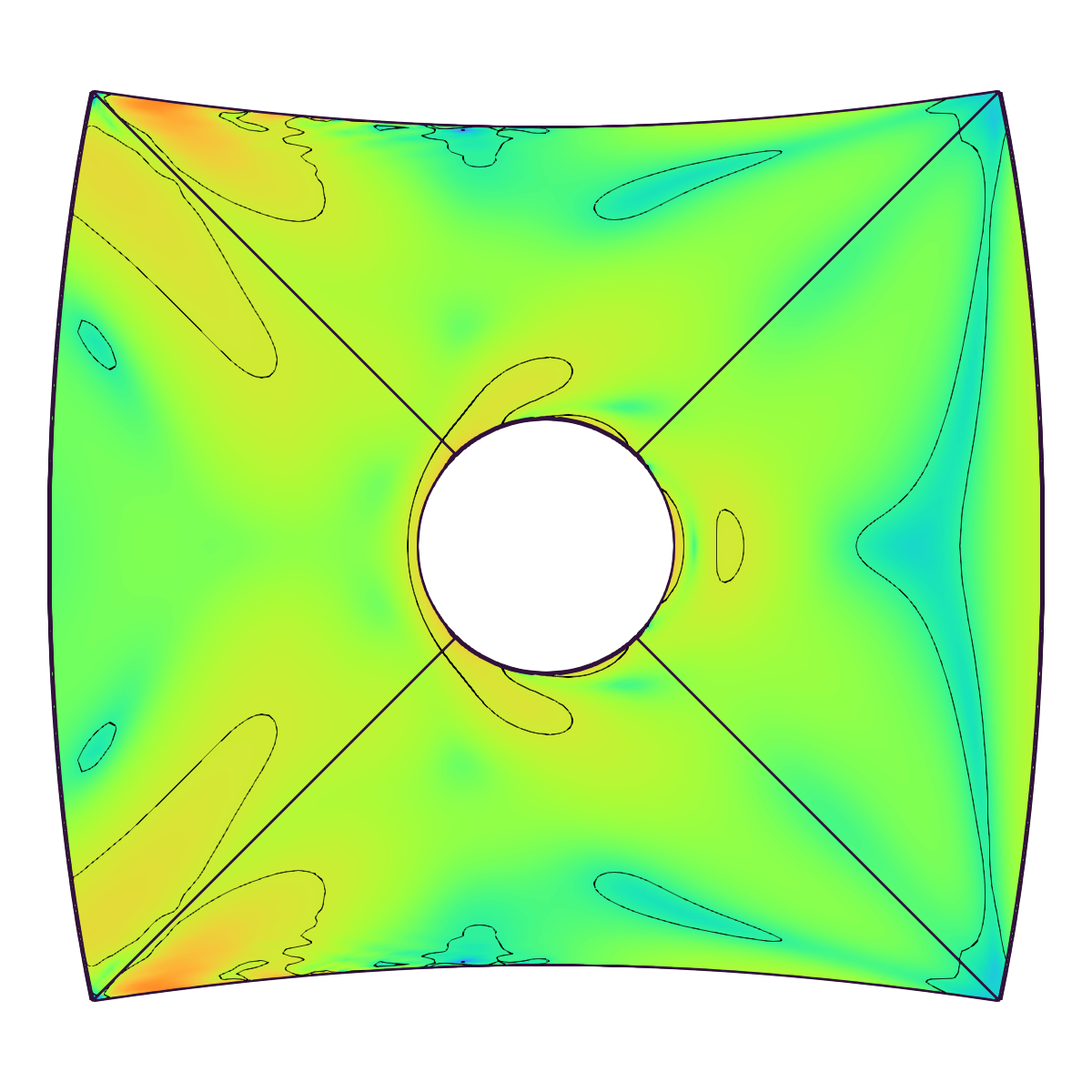}
        \caption{AS-$G^1$ method.}
    \end{subfigure}
    \hfill
    \begin{subfigure}[t]{0.45\linewidth}
        \centering
            \includegraphics[width=0.8\linewidth]{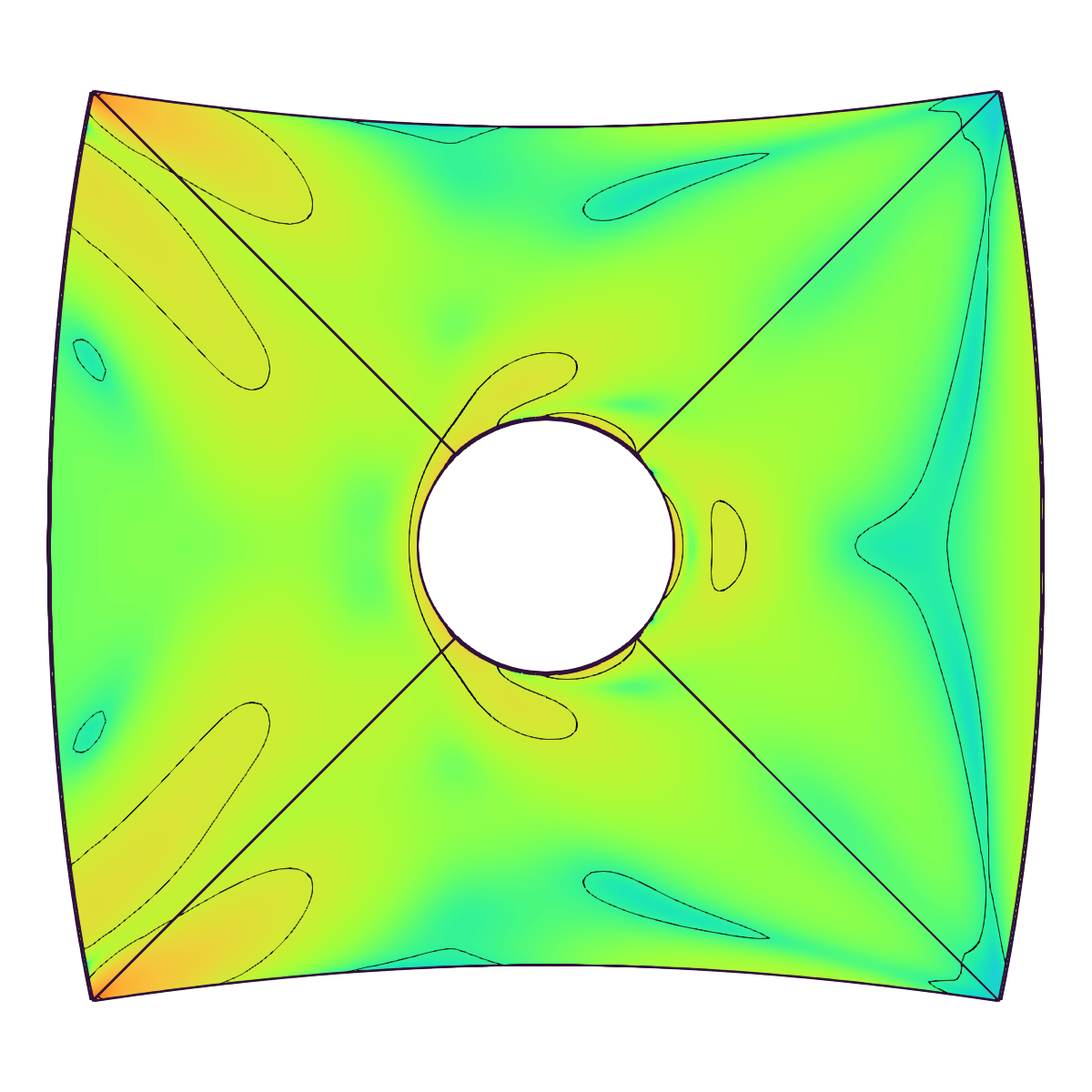}
        \caption{Penalty method.}
    \end{subfigure}
    \caption{Von Mises membrane stress field for the untrimmed hyperboloid geometry with a hole in \cref{fig:hyperboloid} on a $16\times16$ element mesh per patch. The stress field is plotted on the same colour scale as in \cref{fig:hyperboloid_stress}. The contour lines are provided for values $10$, $10^2$, $10^3$, $10^4$ and $10^5\:[\text{MPa}]$.} 
    \label{fig:hyperboloid_hole_stress}
\end{figure}

We perform on the resulting hyperboloid shell with hole numerical tests as in the previous example in \cref{subsec:hyperboloid} by selecting the same geometric and material properties for the shell, i.e. the length by $L=1\:[\text{m}]$, the thickness by $t=0.01\:[\text{m}]$, the Young's modulus by $E=2 \cdot 10^{11}\:[\text{N}/\text{m}^2]$ and the Poisson's ratio by $\nu = 0.3\:[-]$, again by fully clamping the shell on the west side at $x_1 = -L/2$ by and applying a uniformly distributed vertical load~$p_{x_3} = -8000 \cdot t \:[\text{N}/\text{m}^2]$. \cref{fig:hyperboloid_hole_results} presents the results including a convergence analysis for studying the vertical displacement at the reference point $A=(x^1,x^2,x^3) = (0,0,0)$ for degrees~$p=4,5$ and maximal possible regularity $\reg=p-2$ of the spline space~$\mathcal{A}$ by performing $h$-refinement. For this example no reference solution is provided, since no single-patch equivalent for this case exists, hence a penalty coupling \cite{Herrema2019} with parameter $\alpha=10^2$ is used as reference, which appeared to be the optimal choice for displacement and stress reconstruction. Again, we see that the solutions for our AS-$G^1$ method converge, but with a lower rate compared to the penalty method. The speed of convergence grows by increasing the degree $p$. In addition, the stress field for the untrimmed hyperboloid with a hole by using our AS-$G^1$ method and the penalty method are presented in \cref{fig:hyperboloid_hole_stress} and show good correspondence with again small artifacts  on the top and bottom boundaries for the stress field from the AS-$G^1$ method. As mentioned before, this could be simply solved in future by using for each boundary edge e.g. the entire space instead of the subspace selected in~\cite{FaJuKaTa22}. The stress field resulting from the penalty method show a good stress field as well, but depends, as mentioned previously, on the choice of the penalty parameter.

\subsection{Post-buckling of an L-shaped domain}\label{subsec:L-shape}

\begin{figure}
    \centering
\begin{minipage}[t]{0.45\linewidth}
    \centering
    \resizebox{\linewidth}{!}{\includegraphics{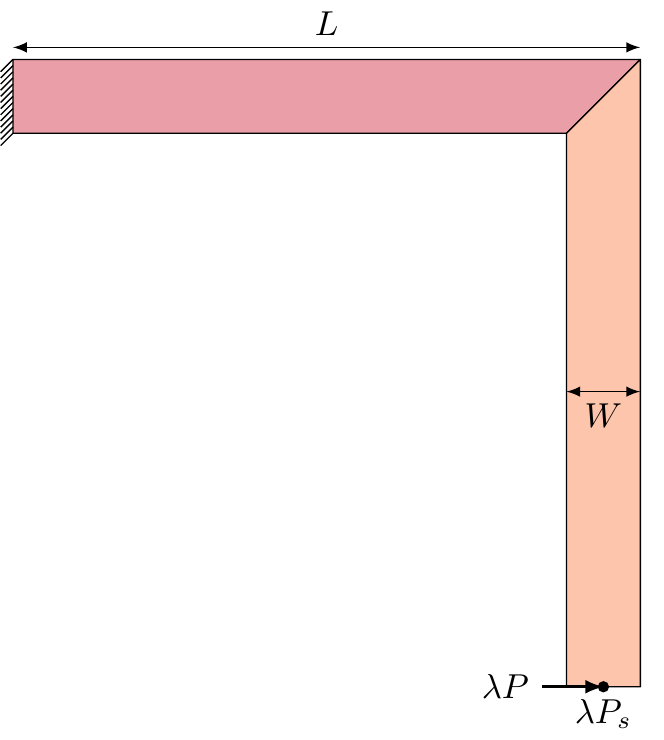}}
    \caption{Geometry of the L-shaped domain with length $L=255\:[\text{mm}]$ and width $W=30\:[\text{mm}]$ consisting of two bi-linear patches. The left side of the domain is fixed in all directions and rotations and on the bottom right a point load with magnitude $\lambda P$ is applied in horizontal direction. To replicate the example from \cite{Simo1986}, an unbalancing load $\lambda P_s$ in out-of-plane direction is applied to force buckling. The factor $\lambda$ is the load magnification factor used in the arc-length methods to solve the problem.
    }
    \label{fig:L-shape_geometry}
\end{minipage}
\hfill
\begin{minipage}[t]{0.45\linewidth}
\centering
\includegraphics[width=\linewidth]{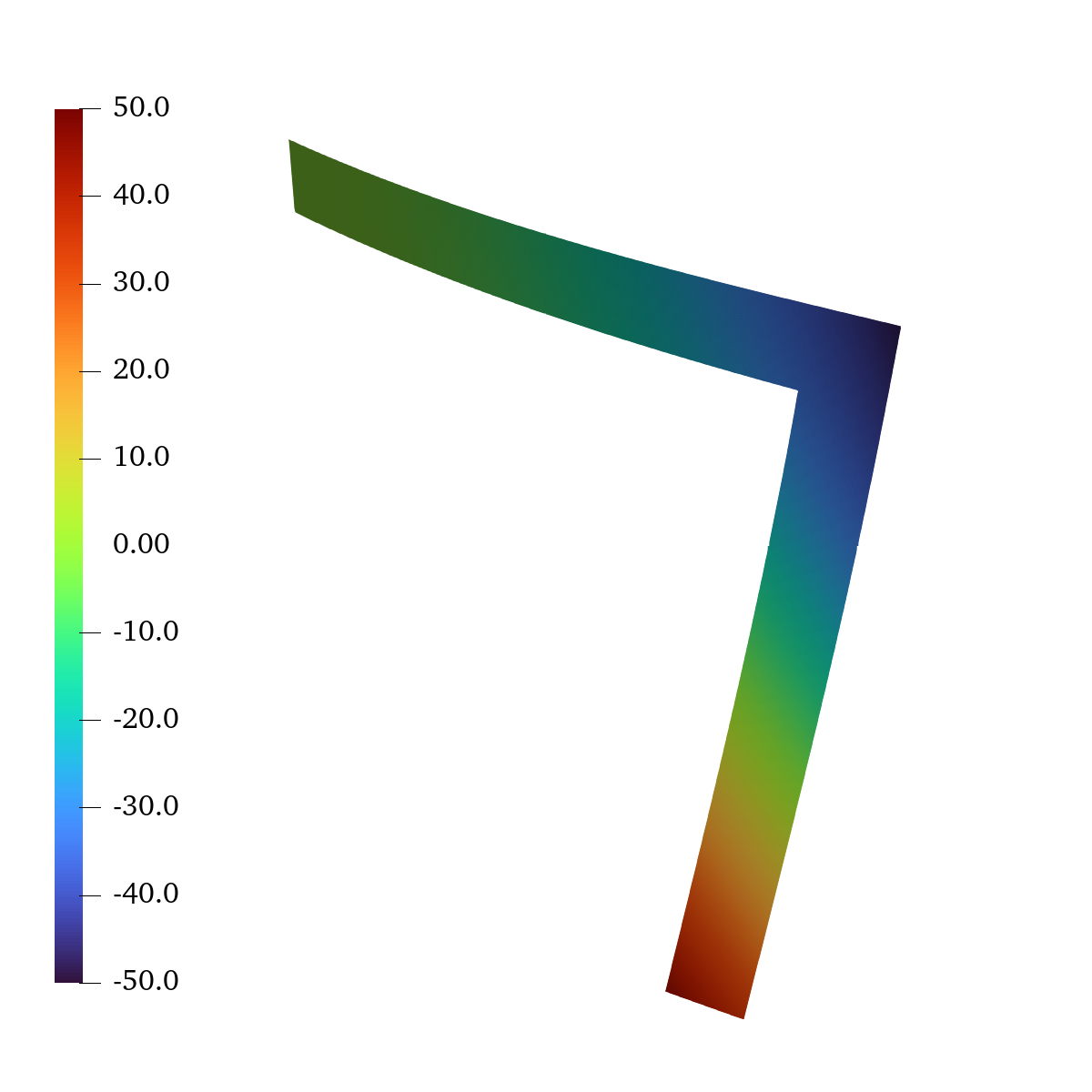}
\caption{Deformed geometry of the L-shaped domain with 2 patches corresponding to \cref{fig:L-shape_geometry} on the last point of the load-displacement curve for the AS-$G^1$ method for $p=5$ in \cref{fig:L-shape_displacements}. The colour scale represents the out-of-plane displacement.}
\label{fig:L-shape_deformed}
\end{minipage}
\end{figure}

We apply our developed method to a geometrically non-linear Kirchhoff-Love shell on the post-buckling response of an L-shaped domain, inspired by an example from \cite{Argyris1979} later performed by \cite{Simo1986} using beam and flat shell elements. The problem set-up is depicted in \cref{fig:L-shape_geometry}. The geometry has length $L=255\:[\text{mm}]$ and width $W=30\:[\text{mm}]$ with a thickness of $0.6\:[\text{mm}]$ and Young's modulus $E=71240\:[\text{N}/\text{mm}^2]$ and Poisson ratio $\nu=0.31\:[-]$. The material law is a linear isotropic Saint-Venant Kirchhoff model. Furthermore, Crisfield's arc-length method is employed \cite{Crisfield1981}. The solutions are computed for degrees $p=3,4,5$ with maximum regularity $\reg=p-2$. Furthermore, the results are calculated by using the penalty method \cite{Herrema2019} with penalty parameter $\alpha=10^3$ and degree $p=3$ and regularity $\reg=2$ as a reference since reported benchmark results have been only provided for beam elements.

\begin{figure}
\centering
    \includegraphics[width=0.95\linewidth]{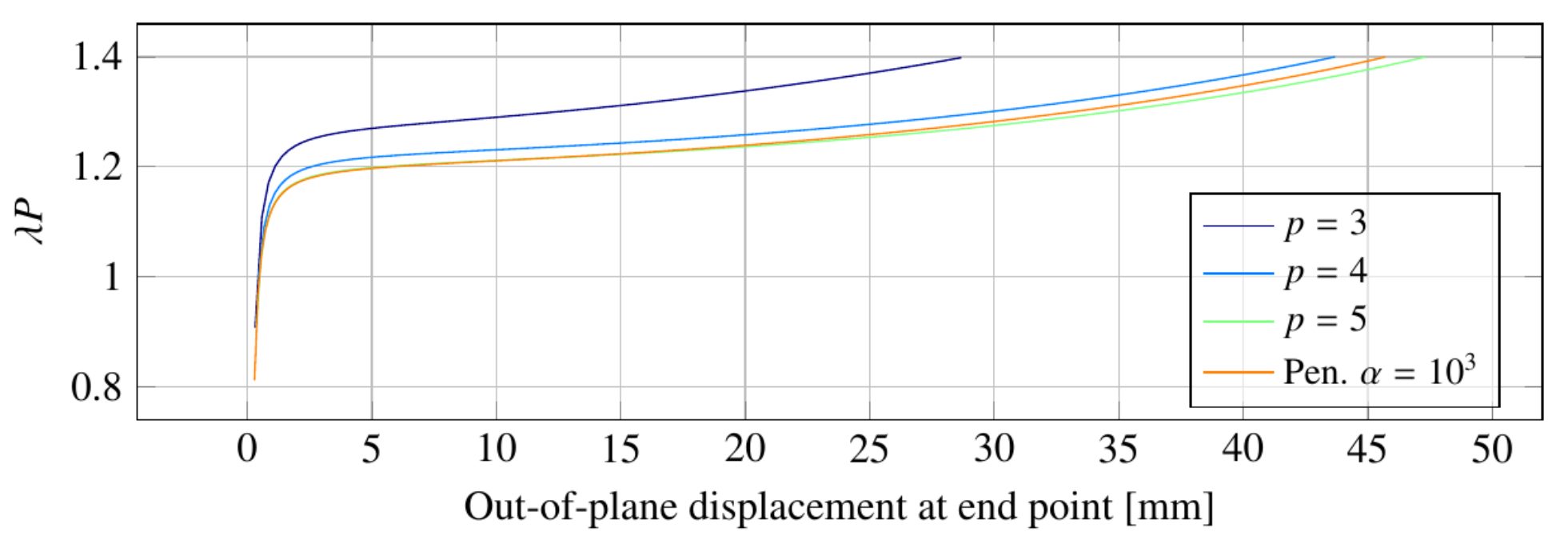}
    \caption{Displacements at the point where the load is applied in \cref{fig:L-shape_geometry}. All results are plotted with the out-of-plane displacement component on the horizontal axis and the load $\lambda P$ on the vertical axis. The penalty parameter used for the penalty method is $\alpha=10^3$.}
    \label{fig:L-shape_displacements}
\end{figure}

The results for the (post-)buckling analysis of the 2 patch L-shaped domain are presented in \cref{fig:L-shape_deformed} and \cref{fig:L-shape_displacements}. As can be seen in \cref{fig:L-shape_displacements}, buckling occurs around $\lambda P=1.19\:[\text{N}]$, which is higher than the value of $\lambda P=1.1453\:[\text{N}]$ as reported by \cite{Argyris1979}, but a converged value for the AS-$G^1$ multi-patch and the penalty coupled multi-patch technique has been tested here. In addition, it can be seen that the load-displacement curves resulting from the AS-$G^1$ method for degree $p=4$ and $p=5$ approach the penalty curve with degree $p=3$. The AS-$G^1$ curve for degree $p=3$, on the other hand, is considered inaccurate.

\subsection{Post-buckling of an L-shaped domain with holes}\label{subsec:L-shape_hole}

\begin{figure}
\centering
\begin{minipage}[t]{0.45\linewidth}
    \centering
    \resizebox{\linewidth}{!}{\includegraphics{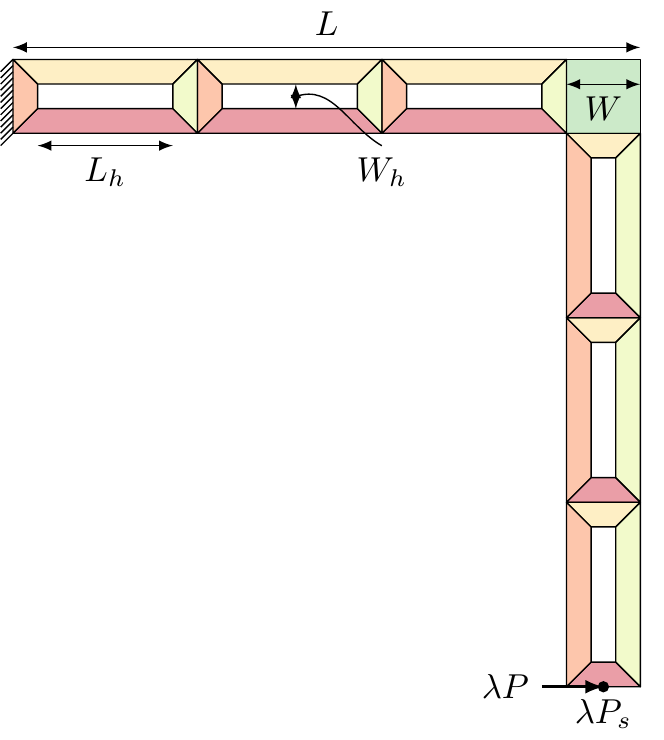}}
    \caption{Geometry of the L-shaped domain with length $L=255\:[\text{mm}]$ and width $W=30\:[\text{mm}]$ and with rectangular holes of size $L_h=55\:[\text{mm}]$ and $W_h=10\:[\text{mm}]$. The geometry consists of 25 bi-linear patches. The setting of the problem is like in 
    \cref{fig:L-shape_geometry}.}
    \label{fig:L-shape_holes_geometry}
\end{minipage}
\hfill
\begin{minipage}[t]{0.45\linewidth}
\centering
\includegraphics[width=\linewidth]{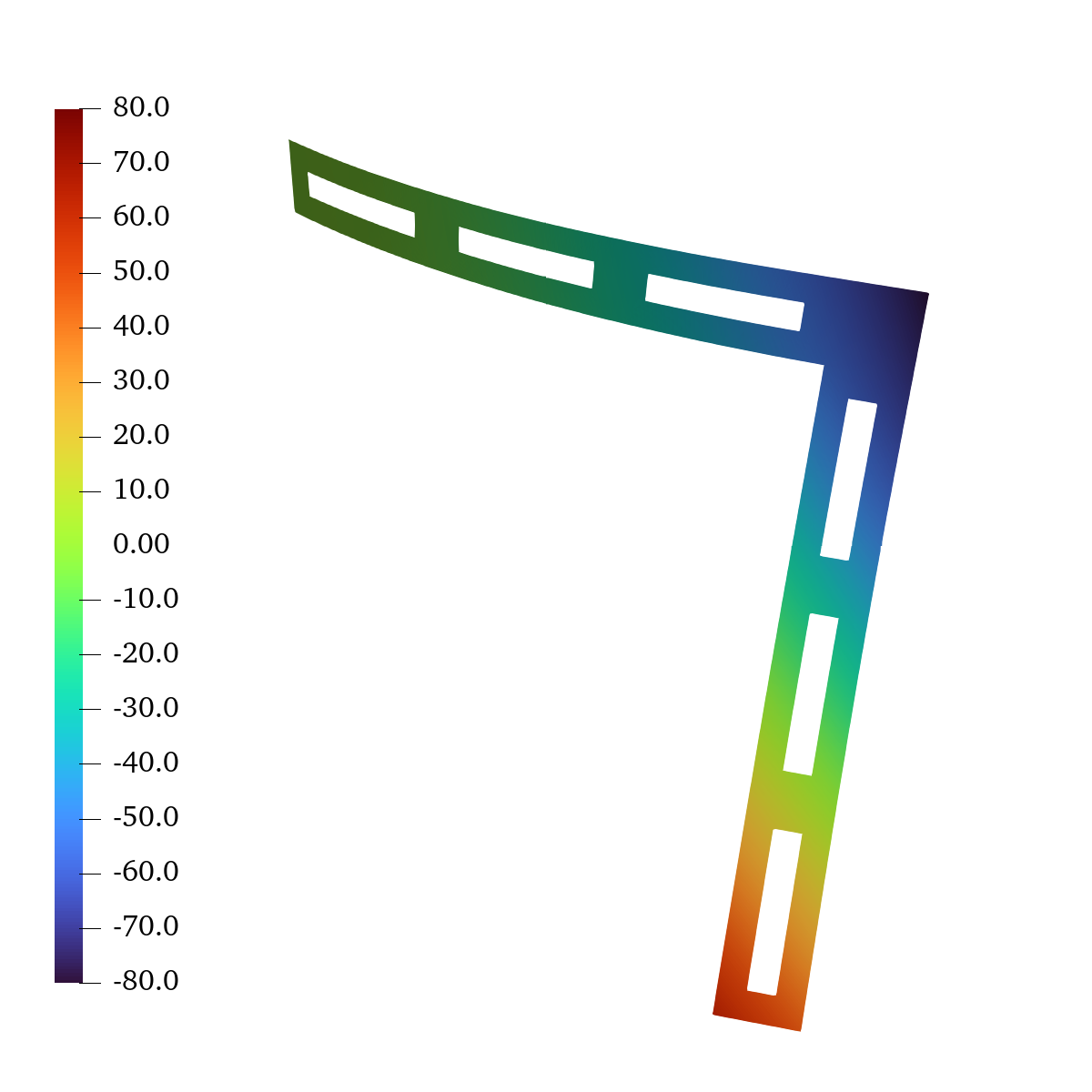}
\caption{Deformed geometry of the L-shaped domain with 25 patches corresponding to \cref{fig:L-shape_holes_geometry} on the last point of the load-displacement curve for the AS-$G^1$ method for $p=4$ in \cref{fig:L-shape_holes_displacements}. The colour scale represents the out-of-plane displacement.}
\label{fig:L-shape_holes_deformed}
\end{minipage}
\end{figure}

As last example, we consider the L-shaped geometry from \cref{subsec:L-shape} with rectangular holes, c.f. \cref{fig:L-shape_holes_geometry}. The material and geometric parameters are the same as in the previous example, additionally with the size of the rectangular holes given by $L_h=55\:[\text{mm}]$ and $W_h=10\:[\text{mm}]$.

\begin{figure}
    \centering
    \includegraphics[width=0.95\linewidth]{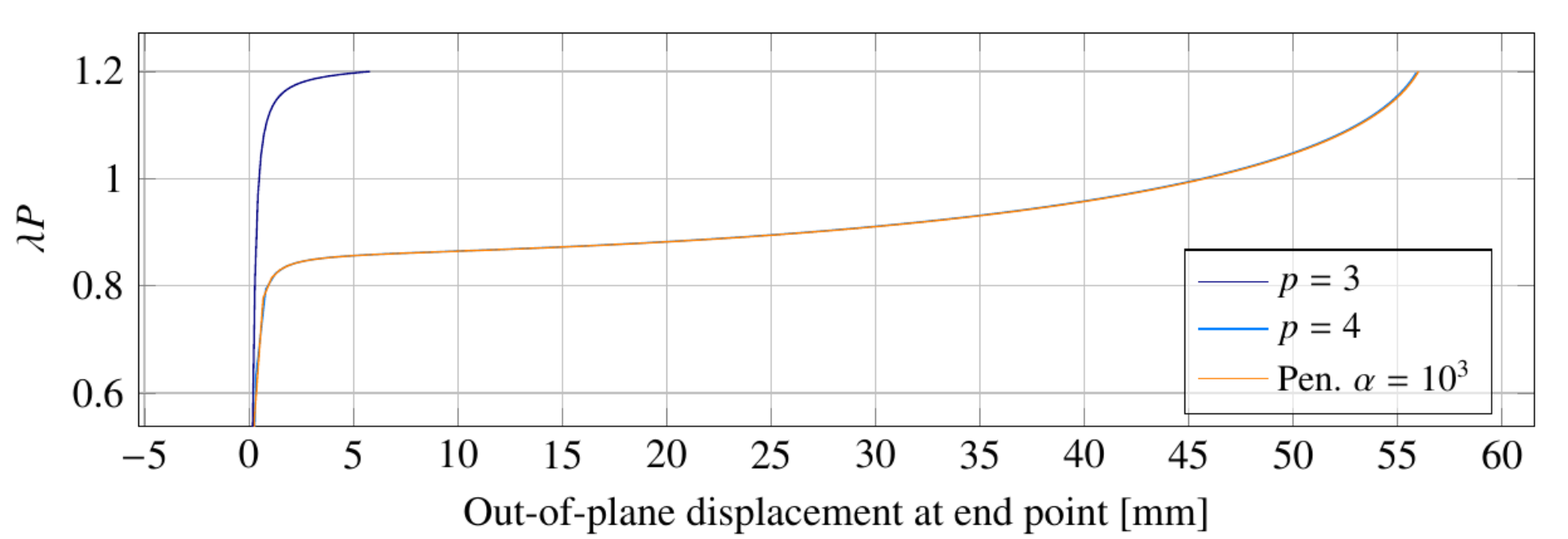}
    \caption{Displacements at the point where the load is applied in \cref{fig:L-shape_holes_geometry}. All results are plotted with the out-of-plane displacement component on the horizontal axis and the load $\lambda P$ on the vertical axis. The penalty parameter used for the penalty method is $\alpha=10^3$.}
    \label{fig:L-shape_holes_displacements}
\end{figure}

The load-displacement curves are given in \cref{fig:L-shape_holes_displacements} for $p=3,4$ and maximum regularity $\reg=p-2$. The curve for $p=5$ is omitted because it overlaps with the $p=4$ curve. Furthermore, the example has been computed by using the penalty method \cite{Herrema2019} with the penalty parameter $\alpha=10^3$. The results show that the AS-$G^1$ technique provides a very good estimate of the load-displacement curve compared to the penalty method for degree $p=4$. For degree $p=3$, however, the results show a large deviation with the other lines, similarly to the results in the previous section.

\section{Conclusion} \label{sec:conclusion}

We presented an isogeometric method for the analysis of complex Kirchhoff-Love shells, where the mid-surface of the shells is approximated by a particular class of $G^1$-smooth multi-patch surfaces, called AS-$G^1$~\cite{CoSaTa16}. This class of multi-patch surfaces allows to use for the discretisation space of the considered shell the globally $C^1$-smooth isogeometric multi-patch spline space with optimal polynomial reproduction properties~\cite{FaJuKaTa22}, whose construction is simple and uniform. Our proposed method further relies on the usage of the Kirchhoff-Love shell formulation~\cite{kiendl-bletzinger-linhard-09}, which can be directly applied due to the globally $C^1$-smoothness of the employed discretisation space. 

We used the developed isogeometric technique to study several linear and non-linear benchmark problems for geometrically complex multi-patch structures. We compared our results to single patch solutions (where possible) and to results obtained by the penalty method \cite{Herrema2019}. We get good results with our AS-$G^1$ approach for displacements and stresses. The convergence of our method is slower than for the penalty method, which is a disadvantage, but the advantages of our AS-$G^1$ approach are that strong $C^1$ coupling is guaranteed, there is no need for the selection of a penalty parameter and that geometries can be handled which are not possible by some other methods like \cite{Toshniwal2017}.

The slower convergence of our method compared to the single patch case and to the penalty method is caused by the fact that the used $C^1$-smooth space~\cite{FaJuKaTa22} is a subspace of the complex, entire $C^1$-smooth space over the considered AS-$G^1$ multi-patch mid-surface. A first possible future work could be now to enlarge the space~\cite{FaJuKaTa22} by constructing the entire $C^1$-smooth space or by partly increasing the degree of the space along the edges, which should considerably increase the speed of convergence. Since the employed $C^1$-smooth spline space~\cite{FaJuKaTa22} allows the construction of a sequence of nested spaces, one further possible future work could be the extension of our approach to an adaptive method to perform local refinement for the analysis of multi-patch Kirchhoff-Love shells. This would allow a reduction of the needed degrees of freedom for solving the problems with an approximation error of similar magnitude.

\section{Acknowledgements} \label{sec:acknowledgements}

The authors wish to thank the anonymous reviewers for their comments that helped to improve the paper. A. Farahat and M. Kapl have been supported by the Austrian Science Fund (FWF) through the project P~33023-N. H.M. Verhelst is grateful for the funding from Delft University of Technology. 
J. Kiendl has received funding from the European Research Council (ERC) under the European Union’s Horizon 2020 research and innovation program (grant agreement No 864482). 
Additionally, the authors are grateful for the support from the developers of the Geometry + Simulation Modules, in particular from A. Mantzaflaris 
(Inria Sophia Antipolis-M\'editerran\'ee).

\section{CRediT Authorship Contributions} \label{sec:credit}
\textbf{Andrea Farahat}: Conceptualization (lead), Formal Analysis (equal), Investigation (equal), Methodology (lead), Software, Validation (equal), Visualization (equal), Writing – Original Draft Preparation (lead), Writing – Review \& Editing (supporting). \textbf{Hugo M. Verhelst}: Conceptualization, Formal Analysis (equal), Investigation (equal), Methodology, Software (lead), Validation (equal), Visualization (equal), Writing – Original Draft Preparation, Writing – Review \& Editing
(lead). \textbf{Josef Kiendl}: Conceptualization, Investigation (supporting), Methodology (supporting), Supervision (equal), Writing – Original Draft Preparation, Writing – Review \& Editing. \textbf{Mario Kapl}: Conceptualization, Funding Acquisition, Investigation (supporting), Methodology (supporting), Project Administration (lead), Supervision (equal), Writing – Original Draft Preparation, Writing – Review \& Editing.

\appendix

\section{Construction of the basis of the $C^1$-smooth subspace~$\mathcal{A}$}  \label{app:construction_space_A}
We will summarize the construction of the basis functions of the $C^1$-smooth isogeometric spline space~$\mathcal{A}$ introduced in~\cite{FaJuKaTa22}, and will further give some implementation details for the functions. Before, some required preliminaries will be presented. 

\subsection{Local parameterisations in standard form} \label{appsec:standardform}

For each edge~$\Sigma^{(i)}$, $i \in \mathcal{I}_{\Sigma}$, or vertex~$\vb{x}^{(i)}$, $i \in \mathcal{I}_{\chi}$, the patch parameterisations~$\pp^{(i_j)}$ in the neighborhood of the edge or vertex can be always locally (re)parameterised (if needed) into the so-called \emph{standard form}~\cite{FaJuKaTa22}. 

\begin{figure}[htbp]
\centering
   \begin{subfigure}[t]{0.43\linewidth}
    \centering
    \includegraphics[width=\linewidth]{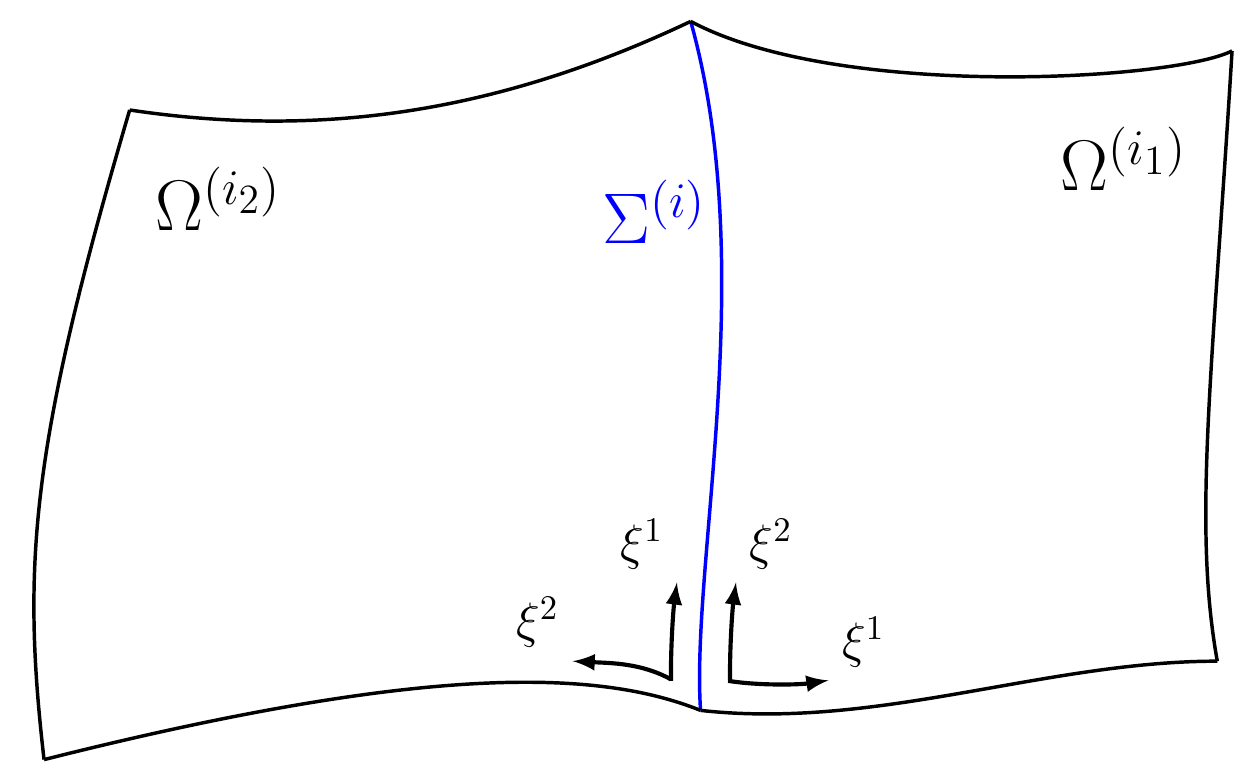}
\end{subfigure}
\begin{subfigure}[t]{0.43\linewidth}
    \centering
    \includegraphics[width=\linewidth]{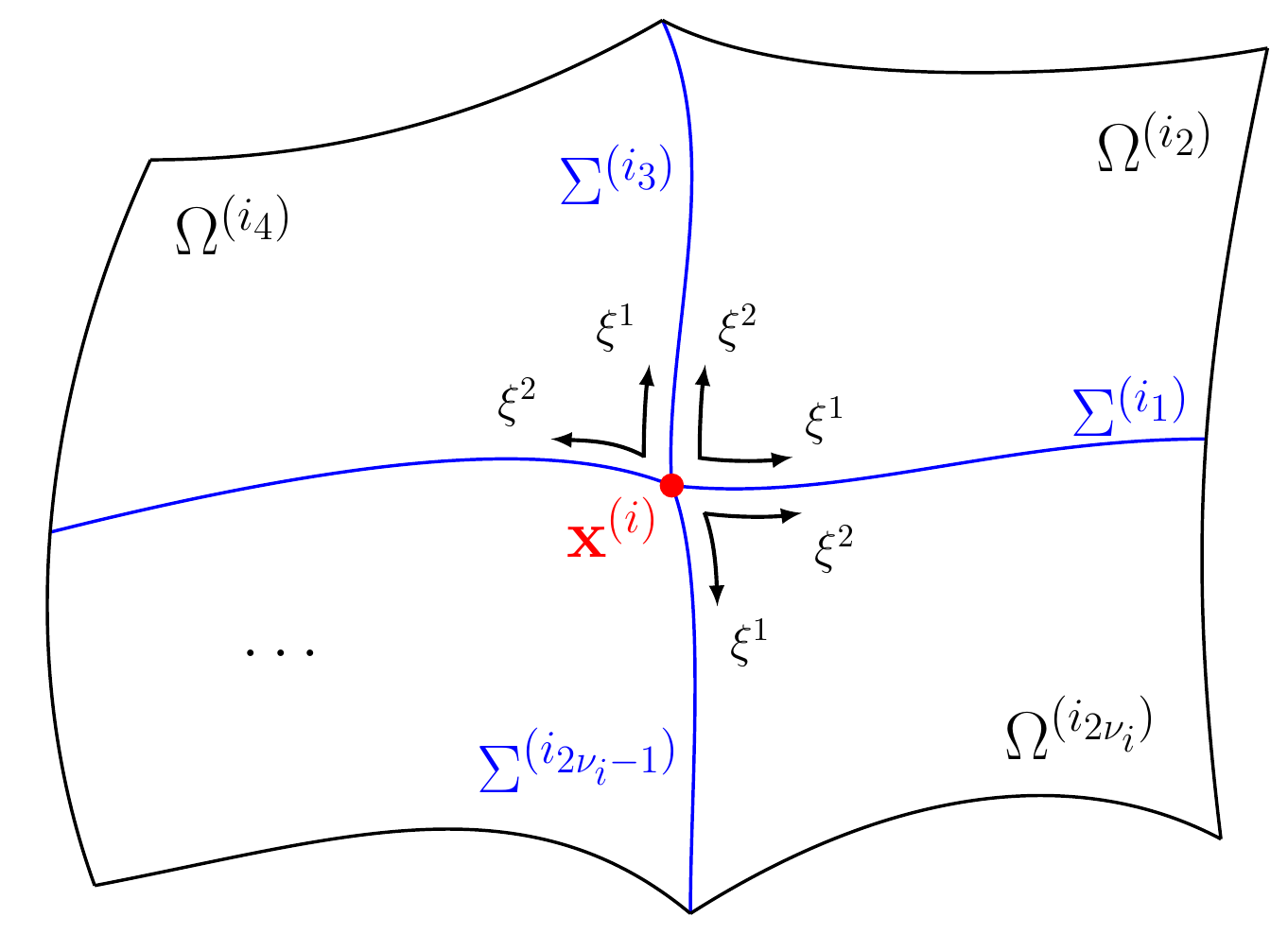}
\end{subfigure}
\caption{Parameterisation in standard form in case of an interface curve (left) and of an inner vertex (right)}
\label{fig:stdParam}
\end{figure}

\paragraph{Standard form for an edge} 
In case of an interface curve $\Sigma^{(i)}$, $i \in \mathcal{I}_{\Sigma}^{\circ}$, with $\Sigma^{(i)} \subset \overline{\Omega^{(i_1)}} \cap \overline{\Omega^{(i_2)}}$, $i_1,i_2 \in \mathcal{I}_{\Omega}$, the patch parameterisations~$\pp^{(i_1)}$ and $\pp^{(i_2)}$ are (re)parameterised in such a way that the curve $\Sigma^{(i)}$ is given by
\[
\pp^{(i_1)}(0,\xi) = \pp^{(i_2)}(\xi,0), \quad \xi \in (0,1), 
\]
see \cref{fig:stdParam} (left).
Similarly, in case of a boundary curve $\Sigma^{(i)}  \subset \overline{\Omega^{(i_1)}}$, $i \in \mathcal{I}_{\Sigma}^{\Gamma}$, the curve $\Sigma^{(i)}$ is just taken as
\[
\Sigma^{(i)}= \{ \pp^{(i_1)}(0,\xi) \mbox{}: \mbox{}\xi \in (0,1) \}. 
\]

\paragraph{Standard form for a vertex}
In case of an inner vertex~$\vb{x}^{(i)}$, $i \in \mathcal{I}_{\chi}^{\circ}$, with patch valence~$\nu_i$, assuming that the patches and interface curves around the vertex~$\vb{x}^{(i)}$ are labeled in counterclockwise order as $\Sigma^{(i_1)}$, $\Omega^{(i_2)}$, $\Sigma^{(i_3)}$, $\ldots$, $\Sigma^{(i_{2{\nu_i-1}})}$, $\Omega^{(i_{2 \nu_i})}$ (and further setting $\Sigma^{(i_{2\nu_i +1})}=\Sigma^{(i_1)}$ and $\Omega^{(i_0)} = \Omega^{(i_{2 \nu_i})}$), the associated patch parameterisations~$\pp^{(i_{2 \ell})}$, $\ell =1,\ldots, \nu_i $, are (re)parameterized in such a way that the interface curve $\Sigma^{(i_1)}$ is given by
\[
\pp^{(i_{2 \nu_i})} (0,\xi) = \pp^{(i_{2})} (\xi,0), \mbox{ } \xi \in (0,1),
\]
and the interface curves~$\Sigma^{(i_{2 \ell +1})}$, $\ell=1, \ldots, \nu_i -1$, are given by
\[
\pp^{(i_{2 \ell})} (0,\xi) = \pp^{(i_{2 \ell +2})} (\xi,0), \mbox{ } \xi \in (0,1),
\]
which leads to
\[
\vb{x}^{(i)} = \pp^{(i_{2})} (0,0) = \ldots = \pp^{(i_{2 \nu_i})} (0,0),
\]
cf. \cref{fig:stdParam} (right). In addition, the patch parameterisations~$\pp^{(i_{2 \ell})}$, $\ell =1,\ldots, \nu_i $, are collectively rotated in such a way that their normal vectors at the vertex~$\vb{x}^{(i)}$ are parallel to the $x_3$-axis. Similarly, the patch parameterisations~$\pp^{(i_{2\ell})}$, $\ell=1,\ldots, \nu_i$, around a boundary vertex~$\vb{x}^{(i)}$, $i \in \mathcal{I}_{\Sigma}^{\Gamma}$, with patch valence~$\nu_i$ are (re)parameterised and rotated, where the curve~$\Sigma^{(i_1)}$ and the additional curve~$\Sigma^{(i_{2 \nu_{i}+1})}$ are the two boundary curves.

\subsection{Construction of the basis functions} \label{appsec:construction}

We present the construction of the patch, edge and vertex functions for the single patches, edges and vertices, respectively. This will require the use of the functions $M_i^{p,\reg}$, $i=0,1$, $M_i^{p,\reg+1}$, $i=0,1,2$, and $M_{i}^{p-1,\reg}$, $i=0,1$, which are defined as
\[
 M_0^{p,\reg}(\xi)=\sum_{j=0}^{1} N_j^{p,\reg}(\xi), \mbox{ } M_1^{p,\reg}(\xi)=\frac{h}{p}N_1^{p,\reg}(\xi),
\]
\[
 M_0^{p-1,\reg}(\xi)=\sum_{j=0}^{1} N_j^{p-1,\reg}(\xi), \mbox{ } M_1^{p-1,\reg}(\xi)=\frac{h}{p}N_1^{p-1,\reg}(\xi),
\]
and
\[
 M_0^{p,\reg+1}(\xi)=\sum_{j=0}^{2} N_j^{p,\reg+1}(\xi), \mbox{ } M_1^{p,\reg+1}(\xi)=\frac{h}{p}\sum_{j=1}^2\vartheta(j) N_j^{p,\reg+1}(\xi), \mbox{ } 
 M_2^{p,\reg+1}(\xi)=\frac{h^2 \mu}{p(p-1)}N_2^{p,\reg+1}(\xi)
\]
with $\vartheta(j)=j$ and $\mu=1$ for $\reg < p-2$, and $\vartheta(j)=2j-1$ and $\mu=2$ for $\reg = p-2$.

\paragraph{Patch functions}
For each patch~$\Omega^{(i)}$, $i \in \mathcal{I}_{\Omega}$, the patch functions~$\phi_{(j_1,j_2)}^{\Omega^{(i)}}$, $j_1,j_2 \in \{2, \ldots, n-3\}$, are constructed as
\begin{equation*}
\phi_{(j_1,j_2)}^{\Omega^{(i)}}(\vb{x}) = 
\begin{cases}
\left( f_{\Omega^{(i)},(j_1,j_2)}^{(i)} \circ \left(\pp^{(i)}\right)^{-1} \right)(\vb{x}) & \mbox{if } \vb{x} \in \overline{\Omega^{(i)}}, \\
0 & \mbox{otherwise},
\end{cases}
\end{equation*}
with 
\[
f_{\Omega^{(i)},(j_1,j_2)}^{(i)}(\xi^1,\xi^2)=N^{\vb{p},\breg}_{(j_1,j_2)}(\xi^1,\xi^2).
\]

\paragraph{Edge functions}
For each interface curve~$\Sigma^{(i)}$, $i \in \mathcal{I}_{\Sigma}^{\circ}$, locally parameterised in standard form, the edge  functions~$\phi_{(j_1,j_2)}^{\Sigma^{(i)}} $, $j_1 \in \{3-j_2,  \ldots, n_{j_2}-4 +j_2 \}$, $ j_2=0,1$, are defined as
\begin{equation*}
\phi_{(j_1,j_2)}^{\Sigma^{(i)}}(\vb{x}) = 
\begin{cases}
\left( f_{\Sigma^{(i)},(j_1,j_2)}^{(i_1)} \circ \left(\pp^{(i_1)}\right)^{-1} \right)(\vb{x}) & \mbox{if } \vb{x} \in \overline{\Omega^{(i_1)}}, \\
\left( f_{\Sigma^{(i)},(j_1,j_2)}^{(i_2)} \circ \left(\pp^{(i_2)}\right)^{-1} \right)(\vb{x}) & \mbox{if } \vb{x} \in \overline{\Omega^{(i_2)}}, \\
0 & \mbox{otherwise},
\end{cases}
\end{equation*}
with
\begin{equation} \label{eq:f0_edge}
\begin{aligned}
f_{\Sigma^{(i)},(j_1,0)}^{(i_1)}(\xi^1,\xi^2)  =  N_{j_1}^{p,\reg +1}(\xi^2)M_{0}^{p,\reg}(\xi^1) - \beta^{(i,i_1)}(\xi^2) \left( N_{j_1}^{p,\reg +1} \right)'(\xi^2)M_{1}^{p,\reg}(\xi^1), \\
f_{\Sigma^{(i)},(j_1,0)}^{(i_2)}(\xi^1,\xi^2)  =  N_{j_1}^{p,\reg +1}(\xi^1)M_{0}^{p,\reg}(\xi^2) - \beta^{(i,i_2)}(\xi^1) \left( N_{j_1}^{p,\reg +1} \right)'(\xi^1)M_{1}^{p,\reg}(\xi^2),
\end{aligned}
\end{equation}
and
\begin{equation} \label{eq:f1_edge}
\begin{aligned}
f_{\Sigma^{(i)},(j_1,1)}^{(i_1)}(\xi^1,\xi^2)  =  \alpha^{(i,i_1)}(\xi^2)N_{j_1}^{p-1,\reg}(\xi^1), \\
f_{\Sigma^{(i)},(j_1,1)}^{(i_2)}(\xi^1,\xi^2)  = -\alpha^{(i,i_2)}(\xi^1)N_{j_1}^{p-1,\reg}(\xi^2).
\end{aligned}
\end{equation}
In case of a boundary curve~$\Sigma^{(i)}$, $i \in \mathcal{I}_{\Sigma}^{\Gamma}$, locally given in standard form, the edge  functions~$\phi_{(j_1,j_2)}^{\Sigma^{(i)}} $, $j_1 \in \{3-j_2,  \ldots, n_{j_2}-4 +j_2 \}$, $\; j_2=0,1$, are equal to
\begin{equation*}
\phi_{(j_1,j_2)}^{\Sigma^{(i)}}(\vb{x}) = 
\begin{cases}
\left( f_{\Sigma^{(i)},(j_1,j_2)}^{(i_1)} \circ \left(\pp^{(i_1)}\right)^{-1} \right)(\vb{x}) & \mbox{if } \vb{x} \in \overline{\Omega^{(i_1)}}, \\
0 & \mbox{otherwise},
\end{cases}
\end{equation*}
with the functions $f_{\Sigma^{(i)},(j_1,0)}^{(i_1)}$ and $f_{\Sigma^{(i)},(j_1,1)}^{(i_1)}$ as given in \eqref{eq:f0_edge} and \eqref{eq:f1_edge}, respectively. 

\paragraph{Vertex functions}
Let $\vb{x}^{(i)}$, $i \in \mathcal{I}_{\chi}$, be a vertex with patch valence~$\nu_i$ locally given in standard form. We denote by $\pp_{P}^{(i_{2 \ell})}$, $\ell=1,\ldots,\nu_i$, the patch parameterisations $\pp^{(i_{2 \ell})}$ restricted to the first two coordinates. For each interface curve~$\Sigma^{(i_{\ell})}$, we define the vector functions 
\[
\vb{t}^{(i_{\ell})}(\xi) = \partial_2 \pp_P^{(i_{\ell-1})}(0,\xi)=\partial_1 \pp_P^{(i_{\ell+1})}(\xi,0)
\]
and
\begin{equation*}
\begin{array}{lll}
\vb{d}^{(i_{\ell})}(\xi) & = & \frac{1}{\alpha^{(i_\ell,i_{\ell-1})}(\xi)} \left( \partial_1 \pp_P^{(i_{\ell-1})}(0,\xi) + \beta^{(i_{\ell},i_{\ell-1})}(\xi) \partial_2 \pp_P^{(i_{\ell-1})}(0,\xi) \right)  \\
& = & -\frac{1}{\alpha^{(i_\ell,i_{\ell+1})}(\xi)} \left( \partial_2 \pp_P^{(i_{\ell+1})}(\xi,0) + \beta^{(i_{\ell},i_{\ell+1})}(\xi) \partial_1 \pp_P^{(i_{\ell+1})}(\xi,0) \right),
\end{array}
\end{equation*}
and in analogous manner for the possible boundary curves~$\Sigma^{(i_1)}$ and $\Sigma^{(i_{2\nu_i +1})}$. The vertex functions $\phi_{(j_1,j_2)}^{\vb{x}^{(i)}}$, $j_1,j_2=0,1,2$, $j_1+j_2 \leq 2$, are defined as  
\begin{equation*}
\phi_{(j_1,j_2)}^{\vb{x}^{(i)}}(\vb{x}) = 
\begin{cases}
\left( \sigma^{j_1 + j_2} f_{\vb{x}^{(i)},(j_1,j_2)}^{(i_{\ell})} \circ \left(\pp^{(i_{\ell})}\right)^{-1} \right)(\vb{x}) & \mbox{if } \vb{x} \in \overline{\Omega^{(i_{\ell})}}, \mbox{ }\ell=2,4, \ldots, 2 \nu_i, \\
0 & \mbox{otherwise},
\end{cases}
\end{equation*}
with the uniform scaling factor $\sigma = \left (  \frac{h}{ p \, \nu{_i}} \sum_{\ell = 1}^{\nu_i} \|\nabla \pp_{P}^{(i_{2\ell})} (0,0) \| \right ) ^{-1}$. The functions $f_{\vb{x}^{(i)},\vb{j}}^{(i_{\ell})}$, $\vb{j}=(j_1,j_2)$, $j_1,j_2=0,1,2$, $j_1+j_2 \leq 2$, are given as  
\[
f_{\vb{x}^{(i)},\vb{j}}^{(i_{\ell})}(\xi^1,\xi^2)=g_{\vb{x}^{(i)},\vb{j}}^{(i_{\ell-1},i_{\ell})}(\xi^1,\xi^2) + g_{\vb{x}^{(i)},\vb{j}}^{(i_{\ell+1},i_{\ell})}(\xi^1,\xi^2) - g_{\vb{x}^{(i)},\vb{j}}^{(i_{\ell})}(\xi^1,\xi^2)
\]
with the single functions
\begin{eqnarray*}
g_{\vb{x}^{(i)},\vb{j}}^{(i_{\ell+1},i_{\ell})}(\xi^1,\xi^2) & = & 
\sum_{\omega=0}^2 d_{\vb{j},(0,\omega)}^{(i_{\ell+1},i_{\ell})} \left( 
M_{\omega}^{p,\reg +1} (\xi^2)M_{0}^{p,\reg}(\xi^1) - \beta^{(i_{\ell+1},i_{\ell})}(\xi^2) (M_{\omega}^{p,\reg +1 })'(\xi^2) M_{1}^{p,\reg}(\xi^1) \right) \\
& + & \sum_{\omega=0}^1 d_{\vb{j},(1,\omega)}^{(i_{\ell+1},i_{\ell})} \alpha^{(i_{\ell +1},i_{\ell})}(\xi^2) M_{\omega}^{p-1,\reg } (\xi^2)M_{1}^{p,\reg}(\xi^1), \\
g_{\vb{x}^{(i)},\vb{j}}^{(i_{\ell-1},i_{\ell})}(\xi^1,\xi^2) & = & 
\sum_{\omega=0}^2 d_{\vb{j},(0,\omega)}^{(i_{\ell-1},i_{\ell})} \left( 
M_{\omega}^{p,\reg +1} (\xi^1)M_{0}^{p,\reg}(\xi^2) - \beta^{(i_{\ell-1},i_{\ell})}(\xi^1) (M_{\omega}^{p,\reg +1 })'(\xi^1) M_{1}^{p,\reg}(\xi^2) \right) \\
& - & \sum_{\omega=0}^1 d_{\vb{j},(1,\omega)}^{(i_{\ell-1},i_{\ell})} \alpha^{(i_{\ell -1},i_{\ell})}(\xi^1) M_{\omega}^{p-1,\reg } (\xi^1)M_{1}^{p,\reg}(\xi^2),
\end{eqnarray*}
and
\[
g_{\vb{x}^{(i)},\vb{j}}^{(i_{\ell})}(\xi^1,\xi^2) = \sum_{\omega_1=0}^1 \sum_{\omega_2}^1 d_{\vb{j},(\omega_1,\omega_2)}^{(i_{\ell})} M_{\omega_1}^{p, \reg} (\xi^1) M_{\omega_2}^{p, \reg} (\xi^2), 
\]
possessing the coefficients
\[
d_{\vb{j},(0,0)}^{(i_m,i_{\ell})} = \delta_{0 j_1} \delta_{0 j_2}, \mbox{ }
d_{\vb{j},(0,1)}^{(i_m,i_{\ell})} = \vb{b}_{\vb{j}}^{\delta}\,  \vb{t}^{(i_m)}(0), \mbox{ }
d_{\vb{j},(0,2)}^{(i_m,i_{\ell})} = (\vb{t}^{(i_{m})}(0))^T \, H_{\vb{j}}^{\delta} \, \vb{t}^{(i_m)}(0) + \vb{b}_{\vb{j}}^{\delta} \, (\vb{t}^{(i_m)})'(0), 
\]
\[
d_{\vb{j},(1,0)}^{(i_m,i_{\ell})} = \vb{b}_{\vb{j}}^{\delta} \, \vb{d}^{(i_m)}(0), \mbox{ }
d_{\vb{j},(1,1)}^{(i_m,i_{\ell})} = (\vb{t}^{(i_{m})}(0))^T \, H_{\vb{j}}^{\delta} \, \vb{d}^{(i_m)}(0) + \vb{b}_{\vb{j}}^{\delta} \, (\vb{d}^{(i_m)})'(0), 
\]
for $m=\ell-1,\ell+1$, and
\[
d_{\vb{j},(0,0)}^{(i_{\ell})} = \delta_{0 j_1} \delta_{0 j_2}, \mbox{ }
d_{\vb{j},(1,0)}^{(i_{\ell})} = \vb{b}^{\delta}_{\vb{j}} \, \vb{t}^{(i_{\ell-1})}(0), \mbox{ }
d_{\vb{j},(0,1)}^{(i_{\ell})} = \vb{b}^{\delta}_{\vb{j}} \, \vb{t}^{(i_{\ell+1})}(0),
\]
\[
d_{\vb{j},(1,1)}^{(i_{\ell})} = (\vb{t}^{(i_{\ell-1})}(0))^T \, H^{\delta}_{\vb{j}} \, \vb{t}^{(i_{\ell+1})}(0) + \vb{b}^{\delta}_{\vb{j}} \, \partial_1 \partial_2 \pp_P^{(i_{\ell})}(0,0),
\]
with the row vectors $\vb{b}^{\delta}_{\vb{j}} = \begin{bmatrix} \delta_{1 j_1} \delta_{0 j_2} & \delta_{0 j_1} \delta_{1 j_2}  \end{bmatrix}$, and with the matrices
\[
H^{\delta}_{\vb{j}} = \begin{bmatrix}
\delta_{2 j_1} \delta_{0 j_2} & \delta_{1 j_1} \delta_{1 j_2} \\
\delta_{1 j_1} \delta_{1 j_2} & \delta_{0 j_1} \delta_{2 j_2} .
\end{bmatrix} ,
\]
where $\delta_{j_1 j_2}$ is the Kronecker delta.

\subsection{Implementation details}

Due to~$\pp^{(i)} \in (\mathcal{S}_h^{\vb{p}, \breg})^3$, $i \in \mathcal{I}_{\Omega}$, each patch parameterisation~$\pp^{(i)}$ possesses a spline representation of the form
\[
\pp^{(i)}(\xi^1,\xi^2) = \sum_{j_1=0}^{n-1} \sum_{j_2=0}^{n-1} \vb{c}^{(i)}_{(j_1,j_2)} N^{\vb{p},\breg}_{(j_1,j_2)}(\xi^1,\xi^2) =  (\vb{N}^{p,\reg}(\xi^1) )^T \, \vb{C}^{(i)} \, \vb{N}^{p,\reg}(\xi^2),
\]
with the column vector of B-splines $\vb{N}^{p,\reg} = [N_{j}^{p,\reg}]_{j=0,\ldots,n-1}$, and the coefficient matrix~$\vb{C}^{(i)} = [\vb{c}^{(i)}_{(j_1,j_2)}]_{j_1,j_2=0,\ldots,n-1}$, with the single coefficients~$\vb{c}^{(i)}_{(j_1,j_2)} \in \mathbb{R}^3$. Similarly, for each function $\phi \in \mathcal{A}$, the associated spline functions $\phi \circ \pp^{(i)} \in \mathcal{S}_h^{\vb{p}, \breg} $, $i \in \mathcal{I}_{\Omega}$, can be represented as
\begin{equation} \label{eq:representation}
\left( \phi \circ \pp^{(i)} \right)(\xi^1,\xi^2) = \sum_{j_1=0}^{n-1} \sum_{j_2=0}^{n-1} d^{(i)}_{(j_1,j_2)} N^{\vb{p},\breg}_{(j_1,j_2)}(\xi^1,\xi^2) =  (\vb{N}^{p,\reg}(\xi^1) )^T \, D^{(i)} \, \vb{N}^{p,\reg}(\xi^2),
\end{equation}
with the coefficient matrix~$D^{(i)} = [d^{(i)}_{(j_1,j_2)}]_{j_1,j_2=0,\ldots,n-1}$, and the single coefficients~$d^{(i)}_{(j_1,j_2)} \in \mathbb{R}$. Knowing the spline representations~\eqref{eq:representation} for the basis functions of $\mathcal{A}$, we can compute their Bernstein B\'{e}zier representations by means of the so-called B\'{e}zier extraction. This further allow us to implement and integrate the basis functions in isogeometric analysis software such as G+Smo~\cite{JuLaMaMoZu14} or GeoPDEs~\cite{Va16}. Recall that all basis functions are locally supported which means more precisely that in case of a patch function~$\phi_{(j_1,j_2)}^{\Omega^{(i)}}$ just for the patch~$\Omega^{(i)}$, in case of an edge function~$\phi_{(j_1,j_2)}^{\Sigma^{(i)}} $ just for the one patch or the two patches containing the edge~$\Sigma^{(i)}$, and in case of a vertex function~$\phi_{(j_1,j_2)}^{\vb{x}^{(i)}}$ just for the patches containing the vertex~$\vb{x}^{(i)}$, the coefficients in~\eqref{eq:representation} can be non-zero. Below, we will briefly explain how to compute these non-zero coefficients for the single basis functions.

In case of a patch function~$\phi_{(j_1,j_2)}^{\Omega^{(i)}}$, the spline representations~\eqref{eq:representation} are directly obtained, since just one coefficient, namely $d^{(\ell)}_{(j_1,j_2)}=1$ for $\ell=i$, is unequal to zero by definition. In case of an edge function~$\phi_{(j_1,j_2)}^{\Sigma^{(i)}}$ or of a vertex function~$\phi_{(j_1,j_2)}^{\vb{x}^{(i)}}$, the spline representations~\eqref{eq:representation} can be computed by means of the following steps:
\begin{itemize}
    \item At first, we locally (re)parameterise (if needed) the initial patch parameterisations~$\pp^{(i_{\ell})}$ containing the edge~$\Sigma^{(i)}$ or the vertex~$\vb{x}^{(i)}$ into standard form, cf.~\ref{appsec:standardform}. This can be achieved by multiplying the coefficient matrix~$\vb{C}^{(i_{\ell})}$ of each such patch parameterisation~$\pp^{(i_{\ell})}$ by a suitable transformation matrix~$T^{(i_\ell)}$ from the right side, which will accordingly reverse and swap the parametric directions of the patch parameterisations to get the local (re)parameterisation as shown in Fig.~\ref{fig:stdParam}. In case of a vertex function~$\phi_{(j_1,j_2)}^{\vb{x}^{(i)}}$, we additionally multiply each of the three-dimensional coordinates of the transformed coefficient matrix~$\vb{C}^{(i_{\ell})} \, T^{(i_{\ell})}$ by the same suitable rotation matrix from the right side such that the normal vectors of the resulting patch parameterisations at the vertex~$\vb{x}^{(i)}$ are parallel to the $x_3$-axis.
    \item Secondly, we use the obtained patch parameterisations in standard form to construct the edge function~$\phi_{(j_1,j_2)}^{\Sigma^{(i)}}$ or the vertex function~$\phi_{(j_1,j_2)}^{\vb{x}^{(i)}}$ as described in \ref{appsec:construction}.
    \item Thirdly, we compute for each obtained spline function $f_{\Sigma^{(i)},(j_1,j_2)}^{(i_{\ell})} = \phi_{(j_1,j_2)}^{\Sigma^{(i)}} \circ \pp^{(i_{\ell})}$ or $f_{\vb{x}^{(i)},(j_1,j_2)}^{(i_{\ell})} = \phi_{(j_1,j_2)}^{\vb{x}^{(i)}} \circ \pp^{(i_{\ell})}$ the corresponding spline representation~\eqref{eq:representation}, e.g, by means of an interpolation problem.
    \item Finally, we multiply the coefficient matrix~$D^{(i_{\ell})}$ of each of the resulting spline representations~\eqref{eq:representation} of an edge or vertex function from the right side by the inverse of the associated transform matrix~$T^{(i_{\ell})}$, used in the first step, to get the final and desired representation~\eqref{eq:representation}. This will accordingly reverse and swap the parametric directions of the spline functions of an edge or vertex functions in such a way that their parametric directions will coincide again with the ones of the initial patch parameterisations.  
\end{itemize}

\end{document}